\documentclass{amsart}
\usepackage{amsmath,amssymb}
\usepackage[dvips]{graphicx}
\usepackage{subfigure}
\usepackage{psfrag}

\newcommand{\co}{\colon}
\input{macros.cfg}

\newcommand{\entrylabellsmizebul}[1]{{#1}}
\newenvironment{itembul}[1]%
 {\begin{list}{}{\settowidth{\labelwidth}{\entrylabellsmizebul{#1}}%
  \settowidth{\labelsep}{\ }
  \setlength{\topsep}{0pt}
  \setlength{\itemsep}{0pt}
  \setlength{\leftmargin}{\labelwidth}\addtolength{\leftmargin}{\labelsep}%
  }}{\end{list}}

\begin{document}

\title[Involutory Hopf group-coalgebras and flat bundles]{Involutory Hopf group-coalgebras and flat bundles over 3-manifolds}
\author{Alexis Virelizier}
\email{virelizi@math.u-strasbg.fr}
\address{I.R.M.A. - Universit\'e Louis Pasteur - C.N.R.S. \\ 7, rue Ren\'e Descartes \\ 67084 Strasbourg Cedex France}
\begin{abstract}
Given a group $\pi$, we use involutary Hopf \p coalgebras to define a scalar invariant of flat \p bundles over \trois
manifolds. When $\pi=1$, this invariant equals to the one of \trois manifolds constructed by Kuperberg from involutary
Hopf algebras. We give examples which show that this invariant is not trivial.
\end{abstract}
\maketitle

\setcounter{tocdepth}{1} \tableofcontents

\section*{Introduction}

This paper is part of a program recently initiated by Turaev \cite{Tur3}, called \emph{Homotopy quantum field theories}
(HQFT), whose purpose is the study of quantum invariants for maps. A $n$-dimensional HQFT with target a space $X$
consists in associating a vector space $V(g)$ to any map $g\co N \to X$, where $N$ is an oriented closed
$(n-1)$-manifold, and a linear map $L(f) \co  V(f_{\,|\partial_- M}) \to V(f_{\,|\partial_+ M})$ to any map $f\co M \to
X$, where $M$ is an oriented $n$-cobordism with $\partial M= \partial_+ M \cup (-\partial_- M)$. In particular, this
assignment must only depend on the homotopy classes of the maps and must satisfy that to compose two such cobordisms
amounts composing their associated linear maps. When $X$ is reduced to a single point, one recovers the notion of a
topological quantum field theory, as described in~\cite{At1}.

Fix a group $\pi$. A HQFT with target an Eilenberg-Mac Lane space of type $K(\pi,1)$ gives rise to invariants of flat \p
bundles. The ``quantum'' approaches of 3-manifolds invariants can be generalized to this setting to get invariants of
flat \p bundles over \trois manifolds (see \cite{Tur1} for the Reshetikhin-Turaev's one and \cite{Vir3} for the
Hennings-Kauffman-Radford's one). The aim of this paper is to generalize the invariants of \trois manifolds constructed
by Kuperberg \cite{Ku1} from involutory Hopf algebras to invariants of flat \p bundles over \trois manifolds.

The algebraic notion we use to replace that of a Hopf algebra is the notion of a Hopf \p coalgebra, introduced by Turaev
in~\cite{Tur1} and studied by the author in~\cite{Vir2}. Briefly speaking, a Hopf \p coalgebra is a family
$H=\{H_\al\}_{\al \in \pi}$ of algebras (over a field $\Bbbk$) endowed with a comultiplication $\Delta=\{\cp{\al}{\be}\co
H_{\al \be} \to H_\al \otimes H_\be\}_{\al,\be \in \pi}$, a counit $\varepsilon\co  H_1 \to \kk$, and an antipode
$S=\{S_\al\co H_\al \to H_\ali\}_{\al \in \pi}$ which verify some compatibility conditions. The case $\pi=1$ is the
standard setting of Hopf algebras.

Fix a Hopf \p coalgebra $H=\{H_\al\}_{\al \in \pi}$ which is involutory, that is, such that its antipode verifies $S_\ali
S_\al=\id_{H_\al}$ for all $\al \in \pi$. Let $\xi=(p\co \Tilde{M} \to M)$ be a flat \p bundle over a \trois manifold
$M$. Suppose that $\xi$ is pointed and denote its base point by $\tilde{x} \in \Tilde{M}$. Let $x=p(\tilde{x})\in M$.
Then there is a uniquely defined homomorphism $\pi_1(M,x) \to \pi$, called monodromy of $\xi$ at $\tilde{x}$. We define a
scalar invariant $K_H(\xi,\tilde{x})$ of the pointed flat \p bundle $(\xi,\tilde{x})$ as follows: we present the base
space $M$ of $\xi$ by a Heegaard diagram, we color this diagram with $\pi$ by using the monodromy of $\xi$ at
$\tilde{x}$, and we associate to this ``\p colored'' Heegaard diagram some structure constants of $H=\{H_\al\}_{\al \in
\pi}$. The proof of this result consists in showing that the ``\p colored'' Reidemeister-Singer moves report the
equivalence of pointed flat \p bundles over \trois manifolds, and in verifying the invariance of $K_H$ under these moves
by using the properties of involutory Hopf \p coalgebras.

The invariant $K_H$ is not trivial (we give examples of computation).

We study the dependence of $K_H(\xi, \tilde{x})$ in the base point $\tilde{x}$. In particular, we have that $K_H(\xi,
\tilde{x})$ only depends on the path-connected component of $\tilde{x}$ and that $K_H(\xi, \tilde{x})$ is independent of
the choice of $\tilde{x}$ when $\pi$ is abelian.

If $\pi=1$ and $M$ is a \trois manifold, then $K_H(\id_M\co  M \to M)$ coincides with the invariant of $M$ constructed by
Kuperberg \cite{Ku1}.

This paper is organized as follows. In Section~\ref{sect-kup0}, we review properties of Hopf \p coalgebras. In
Section~\ref{sect-kup1}, we construct an invariant of \p colored Heegaard diagrams. In Section~\ref{sect-kup2}, we show
that this invariant yields to an invariant of pointed flat \p bundles over \trois manifolds.
Finally, we give examples in Section~\ref{sect-kup3}.\\

\noindent {\bf Acknowledgements.} The author thanks its advisor Vladimir Turaev for useful suggestions and advices, and
B. Enriquez and L. Vainerman for helpful conversations concerning the algebraic part.\\

Throughout this paper, we let $\pi$ be a group (with neutral element~$1$) and $\Bbbk$ be a field. All algebras are
supposed to be over $\Bbbk$, associative, and unitary. The tensor product $\otimes=\otimes_\Bbbk$ is always assumed to be
over $\Bbbk$. If $U$ and $V$ are \kt spaces, $\sigma_{U,V}\co U \otimes V \to V \otimes U$ will denote the flip defined
by $\sigma_{U,V}(u \otimes v)=v \otimes u$ for all $u \in U$ and $v \in V$.

\section{Hopf group-coalgebras}\label{sect-kup0}
In this section, we review definitions and properties concerning Hopf group-coalgebras. For a detailed treatment, we
refer to \cite{Vir2}.

\subsection{Hopf \protect\p coalgebras}\label{hopfpicoal}
Following~\cite{Tur1}, a \emph{Hopf \p coalgebra} (over $\Bbbk$) is a family $H=\{H_\al\}_{\al \in \pi}$ of \kt algebras
endowed with a family $\Delta=\{\cp{\al}{\be}\co H_{\al \be} \to H_\al \otimes H_\be \}_{\al,\be \in \pi}$ of algebra
homomorphisms (the \emph{comultiplication}), an algebra homomorphism $\varepsilon\co H_1 \to \Bbbk$ (the \emph{counit}),
and a family $S=\{S_\al \co  H_\al \to H_\ali \}_{\al \in \pi}$ of $\Bbbk$-linear maps (the \emph{antipode}) such that,
for all $\al,\be,\ga \in \pi$,
\begin{equation}\label{coass}
(\cp{\al}{\be}\otimes \id_{H_\ga}) \cp{\al \be}{\ga}=(\id_{H_\al} \otimes \cp{\be}{\ga}) \cp{\al}{\be \ga},
\end{equation}
\begin{equation}\label{counit}
(\id_{H_\al} \otimes \varepsilon) \cp{\al}{1}=\id_{H_\al}=(\varepsilon \otimes \id_{H_\al}) \cp{1}{\al},
\end{equation}
and
\begin{equation}\label{antipode}
m_\al (S_\ali \otimes \id_{H_\al}) \cp{\ali}{\al} = \varepsilon \,1_\al =  m_\al (\id_{H_\al} \otimes S_\ali)
\cp{\al}{\ali},
\end{equation}
where $m_\al\co H_\al \otimes H_\al \to H_\al$ and $1_\al \in H_\al$ denote respectively the multiplication and unit
element of $H_\al$.

When $\pi=1$, one recovers the usual notion of a Hopf algebra. In particular $(H_1,m_1,1_1,\cp{1}{1},\varepsilon,S_1)$ is
a Hopf algebra.

Remark that the notion of a Hopf \p coalgebra is not self-dual, and that if $H=\{H_\al\}_{\al \in \pi}$ is a Hopf \p
coalgebra, then $\{ \al \in \pi \, | \, H_\al \neq 0 \}$ is a subgroup of~$\pi$.

A Hopf \p coalgebra $H=\{H_\al\}_{\al \in \pi}$ is said to be of \emph{finite type} if, for all $\al \in \pi$, $H_\al$ is
finite-dimensional (over $\Bbbk$). Note that it does not mean that $\oplus_{\al \in \pi} H_\al$ is finite-dimensional
(unless $H_\al=0$ for all but a finite number of $\al \in \pi$).

The antipode of a Hopf \p coalgebra $H=\{H_\al\}_{\al \in \pi}$ is anti-multiplicative: each $S_\al\co  H_\al \to H_\ali$
is an anti-homomorphism of algebras, and anti-comultiplicative: $\varepsilon S_1=\varepsilon$ and $\cp{\bei}{\ali} S_{\al
\be}=\sigma_{H_\ali,H_\bei}(S_\al \otimes S_\be) \cp{\al}{\be}$ for any $\al,\be \in \pi$, see \cite[Lemma~1.1]{Vir2}.

The antipode $S=\{ S_\al \}_{\al \in \pi}$ of $H=\{H_\al\}_{\al \in \pi}$ is said to be \emph{bijective} if each $S_\al$
is bijective. As for Hopf algebras, the antipode of a finite type Hopf \p coalgebra is always bijective (see
\cite[Corollary~3.7(a)]{Vir2}).

We extend the Sweedler notation for a comultiplication to the setting of a Hopf \p coalgebra $H=\{H_\al\}_{\al \in \pi}$
in the following way: for any $\al,\be \in \pi$ and $h \in H_{\al \be}$, we write $
  \cp{\al}{\be}(h)=\sum_{(h)} \hua \otimes \hdb \in H_\al \otimes H_\be,
$
or shortly, if we leave the summation implicit,
  $
  \cp{\al}{\be}(h)= \hua \otimes \hdb
  $. The coassociativity of $\Delta$ gives that, for any $\al,\be, \ga \in \pi$ and $h \in H_{\al \be \ga}$,
  \begin{equation*}
  h_{(1,\al \be)(1,\al)} \otimes h_{(1,\al \be)(2,\be)} \otimes
  h_{(2,\ga)} = h_{(1,\al)} \otimes h_{(2,\be \ga)(1,\be)} \otimes h_{(2,\be \ga)(2,\ga)}.
  \end{equation*}
This element of $H_\al \otimes H_\be \otimes H_\ga$ is written as $\hua \otimes \hdb \otimes \htg$. By iterating the
procedure, we define inductively $h_{(1,\al_1)} \otimes \dotsb \otimes h_{(n,\al_n)}$ for any $h \in H_{\al_1 \dotsm
\al_n}$.

If $H=\{H_\al\}_{\al \in \pi}$ is a Hopf \p coalgebra with bijective antipode and $H^\opp_\al$ denotes the opposite
algebra to $H_\al$, then $H^\opp=\{H^\opp_\al\}_{\al \in \pi}$, endowed with the comultiplication and counit of $H$ and
with the antipode $S^\opp=\{ S^\opp_\al=S^{-1}_\ali \}_{\al \in \pi}$, is a Hopf \p coalgebra, called \emph{opposite to}
$H$.

Let $H=\{H_\al\}_{\al \in \pi}$ be a Hopf \p coalgebra with bijective antipode. Set $H^\cop_\al=H_\ali$ as an algebra,
$\cpcop{\al}{\be}=\sigma_{H_\bei,H_\ali} \cp{\bei}{\ali}$, $\varepsilon^\cop=\varepsilon$, and $S^\cop_\al=S^{-1}_\al$.
Then $H^\cop=\{H_\al^\cop\}_{\al \in \pi}$ is a Hopf \p coalgebra, called \emph{coopposite to} $H$.

\subsection{The case $\pi$ finite}\label{abstrasctreformulation}
Let us suppose that $\pi$ is a finite group. Recall that the Hopf algebra $\fpi=k^\pi$ of functions on $\pi$ has a basis
$(e_\al\co \pi \to \kk)_{\al \in \pi}$ defined by $e_\al(\be)= \delta_{\al,\be}$ where $\delta_{\al,\al}=1$ and
$\delta_{\al,\be}=0$ if $\al \neq \be$. The structure maps of $\fpi$ are given by $e_\al e_\be  = \delta_{\al,\be} \,
e_\al$, $1_\fpi  = \sum_{\al \in \pi} e_\al$, $\Delta(e_\al) =\sum_{\be \ga=\al} e_\be \otimes e_\ga$,
$\varepsilon(e_\al) =\delta_{\al,1}$, and $S(e_\al) =e_\ali$.

By a \emph{central prolongation} of $\fpi$ we shall mean a Hopf algebra $A$ endowed with a morphism of Hopf algebras
$\fpi \to A$ which sends $\fpi$ into the center of $A$. The morphism $\fpi \to A$ is called the \emph{central map} of
$A$.

Since $\pi$ is finite, any Hopf \p coalgebra $H=\{H_\al\}_{\al \in \pi}$ gives rise to a Hopf algebra
$\tilde{H}=\oplus_{\al \in \pi} H_\al$ with structure maps given by $\tilde{\Delta}|_{H_\al} =\sum_{\be \ga = \al}
\cp{\be}{\ga}$, $\tilde{\varepsilon}|_{H_\al} = \delta_{\al,1}\,\varepsilon$, $\tilde{m}|_{H_\al \otimes H_\be}
=\delta_{\al,\be}\, m_\al$, $\tilde{1} = \sum_{\al \in \pi} 1_\al$, and $\tilde{S} = \sum_{\al \in \pi} S_\al$.

The \kt linear map $\fpi \to \tilde{H}$ defined by $e_\al \mapsto 1_\al$ clearly gives rise to a morphism of Hopf
algebras which sends $\fpi$ into the center of $\tilde{H}$. Hence $\tilde{H}$ is a central prolongation of $\fpi$.

As noticed by Enriquez \cite{Enr1},  the correspondence which assigns to every Hopf \p coalgebra $H=\{H_\al\}_{\al \in
\pi}$ the central prolongation $\tilde{H}$ of $\fpi$ is one-to-one. Indeed, let $(A,m,1,\Delta,\varepsilon,S)$ be a
central prolongation of $\fpi$. Denote by $1_\al \in A$ the image of $e_\al \in \fpi$ under the central map $\fpi \to A$
of $A$. Set $H_\al =A 1_\al$ for any $\al \in \pi$. Since $\fpi \to A$ is a morphism of Hopf algebras and each $1_\al \in
A$ is central, we have that the family $\{H_\al\}_{\al \in \pi}$ is a Hopf \p coalgebra with structure maps given by
$m_\al =1_\al \cdot m|_{H_\al \otimes H_\al}$, $\cp{\al}{\be} = (1_\al \otimes 1_\be) \cdot \Delta|_{H_{\al\be}}$,
$\varepsilon = \varepsilon|_{H_1}$, and $S_\al = 1_\ali \cdot S|_{H_\al}$. Furthermore we have that $\tilde{H}=A$ as a
central prolongation of $\fpi$, where $\tilde{H}=\oplus_{\al \in \pi} H_\al$ is the central prolongation of $\fpi$
associated to $\{H_\al\}_{\al \in \pi}$ as above.

\subsection{\protect\p integrals}\label{s:defiint}
Let us recall that a left (resp.\@ right) integral for a Hopf algebra $(A,\Delta,\varepsilon,S)$ is an element $\Lambda
\in A$ such that $x \Lambda = \varepsilon(x) \Lambda$ (resp.\@ $\Lambda x = \varepsilon(x) \Lambda$) for all $x \in A$. A
left (resp.\@ right) integral for the dual Hopf algebra $A^*$ is a \kt linear form $\lambda\in A^*$ verifying $(\id_A
\otimes \lambda) \Delta(x)=\lambda(x) 1_A$ (resp.\@ $(\lambda \otimes \id_A) \Delta(x)=\lambda(x) 1_A$) for all $x \in
A$.

By a \emph{left} (resp.\@ \emph{right}) \emph{\p integral} for a Hopf \p coalgebra $H=\{H_\al \}_{\al \in \pi}$, we shall
mean a family of \kt linear forms $\lambda=(\la)_{\al \in \pi} \in \Pi_{\al \in \pi} H^*_\al$ such that
\begin{equation}\label{defiint}
(\id_{H_\al} \otimes \lb) \cp{\al}{\be}(x) = \lab(x)  1_\al
       \quad \text{(resp. \;} (\la \otimes \id_{H_\be}) \cp{\al}{\be}(x) = \lab(x) 1_\be\, \text{)}
\end{equation}
for all $\al,\be \in \pi$ and $x \in H_{\al\be}$.

Note that $\lu$ is a usual left (resp.\@ right) integral for the Hopf algebra $H_1^*$.

A \p integral $\lambda=(\la )_{\al \in \pi}$ is said to be \emph{non-zero} if $\lb \neq 0$ for some $\be \in \pi$. Note
that a non-zero \p integral $\lambda=(\la )_{\al \in \pi}$ verifies that $\la \neq 0$ for all $\al \in \pi$ such that
$H_\al \neq 0$ (and in particular $\lu \neq 0$).

It is known that the space of left (resp.\@ right) integrals for a finite-dimensional Hopf algebra is one-dimensional. In
the setting of Hopf \p coalgebras, we also have that the space of left (resp.\@ right) \p integrals for a finite type
Hopf \p coalgebra is one-dimensional (even when $\pi$ is infinite), see \cite[Theorem~3.6]{Vir2}.

\subsection{Semisimplicity} A Hopf \p coalgebra $H=\{H_\al\}_{\al \in \pi}$ is said to
be \emph{semisimple} if each algebra $H_\al$ is semisimple. Note that a necessary condition for $H$ to be semisimple is
that $H_1$ is finite-dimensional (since any infinite-dimensional Hopf algebra over a field is never semisimple,
see~\cite[Corollary 2.7]{sweed2}). When $H$ is of finite type, then $H$ is semisimple if and only if $H_1$ is semisimple
(see \cite[Lemma~5.1]{Vir2}).

\subsection{Cosemisimplicity} The notion of a comodule over a coalgebra may be extended to the setting of Hopf \p coalgebras.
Namely, a \emph{right \p comodule} over a Hopf \p coalgebra $H=\{H_\al\}_{\al \in \pi}$ is a family $M=\{M_\al\}_{\al \in
\pi}$ of \kt spaces endowed with a family $\rho=\{\rh{\al}{\be}\co M_{\al \be} \to M_\al \otimes H_\be \}_{\al, \be \in
\pi}$ of \kt linear maps such that $(\rh{\al}{\be} \otimes \id_{H_\ga}) \rh{\al \be}{\ga}=(\id_{M_\al} \otimes
\cp{\be}{\ga}) \rh{\al}{\be \ga}$ and $(\id_{M_\al} \otimes \varepsilon) \rh{\al}{1}=\id_{M_\al}$ for all $\al,\be,\ga
\in \pi$. A \emph{\p subcomodule} of $M$ is a family $N=\{N_\al\}_{\al \in \pi}$, where $N_\al$ is a \kt subspace of
$M_\al$, such that $\rh{\al}{\be}(N_{\al \be}) \subset N_\al \otimes H_\be$ for all $\al,\be \in \pi$. The notion of sum
and direct sum of a family of \p subcomodules of a right \p comodule may be defined in the obvious way.

A right \p comodule $M=\{M_\al\}_{\al \in \pi}$ is said to be \emph{simple} if it is \emph{non-zero} (i.e., $M_\al \neq
0$ for some $\al \in \pi$) and if it has no \p subcomodules other than $0=\{0\}_{\al \in \pi}$ and itself. A right \p
comodule which is a direct sum of a family of simple \p subcomodules is said to be \emph{cosemisimple}. A Hopf \p
coalgebra is called \emph{cosemisimple} if it is cosemisimple as a right \p comodule over itself (with comultiplication
as structure maps).

There exists (see \cite[Theorem~5.4]{Vir2}) a Hopf \p coalgebra version of the dual Maschke Theorem: a Hopf \p coalgebra
$H=\{H_\al\}_{\al \in \pi}$ is semisimple if and only if there exists a right \p integral $\lambda=( \la )_{\al\in \pi}$
for $H$ such that $\la(1_\al)=1$ for all $\al \in \pi$ with $H_\al \neq 0$. In particular, when $H$ is of finite type, we
have that $H$ is cosemisimple if and only if the Hopf algebra $H_1$ is cosemisimple.

\subsection{Crossed Hopf \protect\p coalgebras}\label{deficro}
The notion of a crossing for a Hopf \p coalgebra is crucial to define the quasitriangularity of a Hopf \p coalgebra
(see~\cite{Tur1,Vir2}). A Hopf \p coalgebra $H=\{H_\al\}_{\al \in \pi}$ is said to be \emph{crossed} if it is endowed
with a family $\varphi=\{\varphi_\be \co  H_\al \to H_{\be \al \bei} \}_{\al,\be \in \pi}$ of algebra isomorphisms (the
\emph{crossing}) such that each $\varphi_\be$ preserves the comultiplication and the counit, i.e., $(\varphi_\be \otimes
\varphi_\be) \cp{\al}{\ga} =\cp{\be \al \bei}{\be \ga \bei} \varphi_\be$ and $\varepsilon \varphi_\be=\varepsilon$ for
all $\al,\be,\ga \in \pi$, and $\varphi$ is \emph{multiplicative} in the sense that $\varphi_{\be \be'}=\varphi_\be
\varphi_{\be'}$ for all $\be, \be' \in \pi$.

One easily verifies that a crossing preserves the antipode, that is, $\varphi_\be S_\al = S_{\be \al \bei} \varphi_\be$
for all $\al,\be \in \pi$.

A particular class of crossed Hopf \p coalgebras is that of Hopf \p coalgebras with $\pi$ abelian: if $\pi$ is an abelian
group, then a Hopf \p coalgebra $H=\{H_\al\}_{\al \in \pi}$ is always crossed (e.g., by taking
$\varphi_\be|_{H_\al}=\id_{H_\al}$).

When $\pi$ is a finite group, the notion of a crossing can be described by using the language of central prolongations of
$\fpi$ (see Section~\ref{abstrasctreformulation}). More precisely, a central prolongation $A$ of $\fpi$ is crossed if it
is endowed with a group homomorphism $\varphi\co  \pi \to \auth(A)$ such that $\pb(1_\al)=1_{\be \al \bei}$ for all
$\al,\be \in \pi$, where $\auth(A)$ is the group of Hopf automorphisms of the Hopf algebra $A$ and $1_\al$ denotes the
image of $e_\al\in \fpi$ under the central map $\fpi \to A$.

\subsection{Involutory Hopf group-coalgebras}\label{sect-involpicoal} In this section we give some results concerning
involutory Hopf \p coalgebras which are used for topological purposes in Sections~\ref{sect-kup1} and \ref{sect-kup2}.

A Hopf \p coalgebra $H=\{H_\al\}_{\al\in \pi}$ is said to be \emph{involutory} if its antipode $S=\{S_\al\}_{\al \in
\pi}$ is such that $S_\ali S_\al=\id_{H_\al}$ for all $\al \in \pi$.

If $A$ is an algebra and $a \in A$, then $r(a) \in \Endo_\Bbbk(A)$ will denote the right multiplication by $a$ defined by
$r(a)(x)=xa$ for any $x \in A$. Moreover, $\Tr$ will denote the usual trace of \kt linear endomorphisms of a \kt space.

\begin{lemma}\label{leminvol}
Let $H=\{H_\al\}_{\al \in \pi}$ be a finite type Hopf $\pi$-coalgebra with antipode $S=\{S_\al\}_{\al \in \pi}$. Let
$\lambda=(\la)_{\al \in \pi}$ be a right $\pi$-integral for $H$ and $\Lambda$ be a left integral for $H_1$ such that
$\lu(\Lambda)=1$. Let $\al \in \pi$. Then
\begin{enumerate}
\renewcommand{\labelenumi}{{\rm (\alph{enumi})}}
 \item $\Tr(f)=\la(S_\ali(\Lambda_{(2,\ali)})f(\Lambda_{(1,\al)}))$ for all $f \in \Endo_\Bbbk(H_\al)$;
 \item $\Tr(r(a) \circ S_\ali S_\al)=\epsilon(\Lambda) \, \la(a)$ for all $a \in H_\al$;
 \item If $H_\al \neq 0$, then $\Tr(S_\ali S_\al) \neq 0$ if and only if $H$ is semisimple and cosemisimple;
 \item If $H_\al \neq 0$, then $\Tr(S_\ali S_\al)=\Tr(S_1^2)$.
\end{enumerate}
\end{lemma}

\begin{proof} To show Part (a), identify $H_\al^* \otimes H_\al$ and $\Endo_\Bbbk(H_\al)$ by $(p \otimes a)(x)=p(x) a$ for all
$p \in H_\al$ and $a,x \in H_\al$. Under this identification, $\Tr(p \otimes a)=p(a)$. Let $f \in \Endo_\Bbbk(H_\al)$. We
may assume that $f=p \otimes a $ for some $p \in H_\al^*$ and $a \in H_\al$. By \cite[Corollary~3.7(b)]{Vir2}, there
exists $b \in H_\al$ such that $p(x)=\la(b x)$ for all $x \in H_\al$. Now $S_\ali^{-1}(b)=\la(b \Lambda_{(1,\al)})
\Lambda_{(2,\ali)}$ by \cite[Lemma~4.3(b)]{Vir2} and so
\begin{equation*}
 b=\la(b \Lambda_{(1,\al)}) S_\ali(\Lambda_{(2,\ali)}) = p(\Lambda_{(1,\al)}) S_\ali(\Lambda_{(2,\ali)}).
\end{equation*}
Therefore
\begin{eqnarray*}
 \Tr(f) & = & p(a) \;\, = \;\, \la(ba)\\
        & = & \la( S_\ali(\Lambda_{(2,\ali)}) p(\Lambda_{(1,\al)}) a) \\
        & = & \la( S_\ali(\Lambda_{(2,\ali)}) f(\Lambda_{(1,\al)})).
\end{eqnarray*}

Let us show Part (b). Let $a \in H_\al$. Then
\begin{eqnarray*}
\Tr(r(a) \circ S_\ali S_\al)
  & = & \la(S_\ali(\Lambda_{(2_\ali)}) S_\ali S_\al (\Lambda_{(1,\al)}) a) \mbox{ \quad by Part (a)} \\
  & = & \la(S_\ali(S_\al (\Lambda_{(1,\al)})\Lambda_{(2_\ali)}) a) \\
  & = & \la(S_\ali(\epsilon(\Lambda) 1_\ali) a) \mbox{ \quad by \eqref{antipode}} \\
  & = & \epsilon(\Lambda) \la(a).
\end{eqnarray*}

To show Part (c), suppose $H_\al \neq 0$. Since $\Tr(S_\ali S_\al)=\epsilon(\Lambda) \la(1_\al)$ (by Part~(b)), one
easily concludes using the facts that $H$ is semisimple if and only if $\epsilon(\Lambda) \neq 0$ (by
\cite[Lemma~5.1]{Vir2} and \cite[Theorem 5.1.8]{sweed}) and $H$ is cosemisimple if and only if $\la(1_\al) \neq 0$ (by
\cite[Theorem~5.4]{Vir2} since $H_\al \neq 0$).

Let us show Part (d). By using \eqref{defiint}, we have $ \lu(1_1) \, 1_\al = (\lu \otimes
\id_{H_\al})\cp{1}{\al}(1_\al)= \la(1_\al) \, 1_\al, $ and so $\la(1_\al)=\lu(1_1)$ since $1_\al \neq 0$ (because $H_\al
\neq 0$). Therefore, by applying Part (b) twice, we obtain that $ \Tr(S_\ali S_\al)=\epsilon(\Lambda) \,
\la(1_\al)=\epsilon(\Lambda) \, \lu(1_1)=\Tr(S_1^2). $
\end{proof}

\begin{lemma}\label{dimeqforinvolcor}
Let $H=\{H_\al\}_{\al \in \pi}$ be a finite type involutory Hopf $\pi$-coalgebra over a field of characteristic $p$. Let
$\al \in \pi$ with $H_\al \neq 0$. If $p=0$ or $p>|\dim H_\al - \dim H_1|$, then $\dim H_\al=\dim H_1$.
\end{lemma}
\begin{proof}
By Lemma~\ref{leminvol}(d), we have $\Tr(S_\ali S_\al)=\Tr(S_1^2)$, that is $(\dim H_\al )1_\kk=(\dim H_1 )1_\kk$ (since
$H$ is involutory). One easily concludes by using the hypothesis on the characteristic of the field $\kk$.
\end{proof}

\begin{lemma}\label{semcosemforinvol}
Let $H=\{H_\al\}_{\al \in \pi}$ be a finite type involutory Hopf $\pi$-coalgebra. Suppose that $\dim H_1 \neq 0$ in the
ground field $\Bbbk$ of $H$. Then $H$ is semisimple and cosemisimple.
\end{lemma}
\begin{proof}
This follows from Lemma~\ref{leminvol}(c), since $\Tr(S_1^2)=\Tr(\id_{H_1})=\dim H_1 \neq 0$.
\end{proof}

\subsection{Diagrammatic formalism of Hopf group-coalgebras}
\begin{figure}[h]
                         \psfrag{m}[][]{$m_\al$}
                         \psfrag{d}[][]{$\Delta_{\alpha,\beta}$}
                         \psfrag{e}[][]{$\varepsilon$}
                         \psfrag{u}[][]{$1_\al$}
                         \psfrag{S}[][]{$S_\alpha$}
                         \scalebox{.9}{\includegraphics{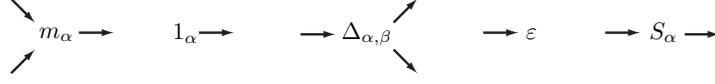}}
     \caption{Symbolic representation of structure maps}
     \label{structuremapdiag}
\end{figure}

The structure maps of a Hopf \p coalgebra $H=\{H_\al\}_{\al \in \pi}$ can be represented symbolically as in \cite{Ku1}.
The multiplications $m_\al\co  H_\al \otimes H_\al \to H_\al$, the units elements $1_\al$, the comultiplication
$\cp{\al}{\be}\co  H_{\al \be} \to H_\al \otimes H_\be$, the counit $\varepsilon\co  H_1 \to \Bbbk$, and the antipode
$S_\al\co H_\al \to H_\ali$ are represented as in Figure~\ref{structuremapdiag}. The inputs (incoming arrows) for the
product symbols are read counterclockwise and the outputs arrows (outgoing arrows) for the comultiplication symbols are
read clockwise.

\begin{figure}[h]
                         \psfrag{m}[][]{$m_\al$}
                         \psfrag{d}[][]{$\Delta_{\alpha_1, \dots, \alpha_n}$}
                         \psfrag{k}[][]{$\Delta_{\alpha_1 \cdots \alpha_{n-1}, \alpha_n}$}
                         \psfrag{l}[][]{$\Delta_{\alpha_1, \alpha_2}$}
                         \psfrag{e}[][]{$\varepsilon$}
                         \psfrag{u}[][]{$1_\al$}
                         \psfrag{a}[][]{$=$}
                         \psfrag{n}[][]{$\Delta_1$}
                         \psfrag{r}[][]{$\Delta_\alpha$}
                         \psfrag{i}[][]{$\id_{H_\alpha}$}
                         \psfrag{T}[][]{$T_\alpha$}
                         \psfrag{C}[][]{$C$}
                         \psfrag{g}[][]{if $\alpha_1 \cdots \alpha_n=1$}
                         \scalebox{.9}{\includegraphics{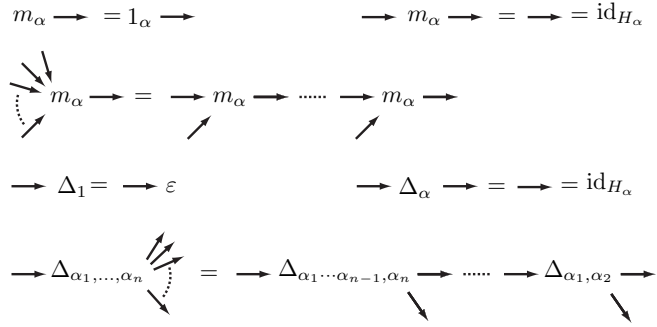}}
     \caption{Diagrammatic abbreviations}
     \label{abbrevdiag}
\end{figure}

In light of the associativity and coassociativity axioms (see Section~\ref{hopfpicoal}), we adopt the abbreviations of
Figure~\ref{abbrevdiag}.

The combinatorics of the diagrams involving such symbolic representations of structure maps may be thought of as (sum of)
products of structure constants. For example, if $(e_i)_i$ is a basis of $H_1$ and $\delta_i^{j,k} \in \Bbbk$ are the
structure constants of $\cp{1}{1}$ defined by $\cp{1}{1}(e_i)=\sum_{j,k} \delta_i^{j,k} e_j \otimes e_k$, then the
element $C \in H_1$ represented in Figure~\ref{diagcandtC} is given by $C=\sum_{i,k} \delta_i^{i,k} e_k$.
\begin{figure}[h]
   \subfigure[$C\in H_1$]{\label{diagcandtC}
                         \psfrag{d}[][]{$\Delta_{1,1}$}
                         \psfrag{C}[][]{$C$}
                         \psfrag{e}[][]{$=$}
                         \scalebox{.9}{\includegraphics{morphC.eps}}} \hspace*{2cm}
   \subfigure[$T_\al\co H_\al \to \kk$]{\label{diagcandtT}
                         \psfrag{m}[][]{$m_\al$}
                         \psfrag{T}[][]{$T_\al = $}
                         \psfrag{e}[][]{$=$}
                         \scalebox{.9}{\includegraphics{morphT.eps}}}
     \caption{}
     \label{diagformalism}
\end{figure}

Similarly, if $(e_i)_i$ is a basis of $H_\al$ and $\mu_{i,j}^k \in \Bbbk$ are the structure constants of $m_\al$ defined
by $m_\al(e_i \otimes e_j)=\sum_k \mu_{i,j}^k e_k$, then the \kt linear form $T_\al\co  H_\al \to \Bbbk$ represented in
Figure~\ref{diagcandtT} is given by $T_\al(e_i)=\sum_k \mu_{k,i}^k$. Note that $T_\al(x)=\Tr(r(x))$ for any $x \in
H_\al$, where $r(x) \in \Endo_\Bbbk(H_\al)$ denotes the right multiplication by $x$ and $\Tr$ is the usual trace of \kt
linear endomorphisms.

Until the end of this section, $H=\{H_\al\}_{\al \in \pi}$ will denote a finite type involutory Hopf \p coalgebra with
$\dim H_1 \neq 0$ in the ground field $\Bbbk$ of $H$.

\begin{lemma}\label{TCtwosided}
$T=(T_\al)_{\al \in \pi}$ is a non-zero two-sided \p integral for $H$ and $C$ is a non-zero two-sided integral for $H_1$
which verify that $T_1(1_1)=\varepsilon(C)=T_1(C)=\dim H_1$. Moreover $S_1(C)=C$ and $T_\ali\circ S_\al=T_\al$ for all
$\al \in \pi$.
\end{lemma}
\begin{proof}
Recall that $H$ is semisimple and cosemisimple (by Lemma~\ref{semcosemforinvol}). Therefore, by \cite[Theorem~5.4 and
Corollary~5.7]{Vir2}, there exists a two-sided \p integral $\lambda=(\la)_{\al \in \pi}$ for $H$ such that $\la(1_\al) =
1$ for all $\al \in \pi$ with $H_\al \neq 0$. Let $\Lambda$ be a left integral for $H_1$ such that $\lu(\Lambda)=1$. By
Lemma~\ref{leminvol}(b), we have that $T_\al(x) = \Tr(r(x))=\varepsilon(\Lambda) \, \la(x)$ for any $x \in H_\al$.
Therefore $T=(T_\al)_{\al \in \pi}$ is a multiple of $\lambda=(\la)_{\al \in \pi}$ and so is a two-sided \p integral for
$H$, which is non-zero since $\varepsilon(\Lambda)\neq 0$ (because $H_1$ is semisimple, see \cite[Theorem 5.1.8]{sweed}).
Likewise $C=\lu(1_1)\, \Lambda=\Lambda$ (by Lemma~\ref{leminvol}(b) applied to the Hopf $1$-coalgebra $H_1^*$) and so $C$
is a non-zero left integral for $H_1$. Moreover $C$ is a right integral for $H_1$ (since $H_1$ is semisimple and so its
integrals are two-sided).

Since $\lu(1_1)=\lu(\Lambda)=1$ and by Lemma~\ref{leminvol}(b), we have that $T_1(C) = T_1(1_1)=\varepsilon(C)
=\varepsilon(\Lambda)=\Tr(\id_{H_1}) =\dim H_1$.

Since $H$ is cosemisimple, \cite[Theorem~4.2(c) and Corollary~5.7]{Vir2} give that $T_\ali \circ S_\al=T_\al$ for all
$\al \in \pi$. Finally, $S_1(C)$ is a left integral for $H_1$ and so there exists $k \in \kk$ such that $S_1(C)=k\, C$.
Now $k=1$ since $C=\Lambda$, $\lu(\Lambda)=1$, and $\lu \circ S_1=\lu$. Hence $S_1(C)=C$.
\end{proof}

\begin{figure}[h]
                         \psfrag{m}[][]{$m_\al$}
                         \psfrag{d}[][]{$\Delta_{\alpha_1, \dots, \alpha_n}$}
                         \psfrag{k}[][]{$\Delta_{\alpha_1 \cdots \alpha_{n-1}, \alpha_n}$}
                         \psfrag{l}[][]{$\Delta_{\alpha_1, \alpha_2}$}
                         \psfrag{e}[][]{$\varepsilon$}
                         \psfrag{u}[][]{$1_\al$}
                         \psfrag{a}[][]{$=$}
                         \psfrag{n}[][]{$\Delta_1$}
                         \psfrag{r}[][]{$\Delta_\alpha$}
                         \psfrag{i}[][]{$\id_{H_\alpha}$}
                         \psfrag{T}[][]{$T_\alpha$}
                         \psfrag{C}[][]{$C$}
                         \psfrag{g}[][]{if $\alpha_1 \cdots \alpha_n=1$}
                         \scalebox{.9}{\includegraphics{abbrevdiag1.eps}}
     \caption{}
     \label{abbrevdiag2}
\end{figure}
\begin{lemma}\label{symabrev}
The two tensors represented by the diagrams of Figure~\ref{abbrevdiag2} are cyclically symmetric.
\end{lemma}
\begin{proof}
Let $\al \in \pi$. Since $(T_\be)_{\be \in \pi}$ is a right \p integral for $H$ (by Lemma~\ref{TCtwosided}) and the Hopf
algebra $H_1$ is semisimple and so unimodular, \cite[Theorem~4.2(a)]{Vir2} gives that
\begin{equation*}
T_\al(xy)=T_\al(S_\ali S_\al(y) x)=T_\al(yx)
\end{equation*}
for all $x,y \in H_\al$. Hence $T_\al(x_1 x_2 \cdots x_n)= T_\al(x_2 \cdots x_n x_1)$ for all $x_1, \dots, x_n \in
H_\al$.

Since $H$ is cosemisimple and $C$ is a left integral for $H_1$, by using \cite[Corollaries~4.4 and 5.7]{Vir2} we have
that, for all $\al \in \pi$,
\begin{equation*}
C_{(1,\al)} \otimes C_{(2,\ali)} = S_\ali S_\al(C_{(2,\al)})1_\al \otimes C_{(1,\ali)} =C_{(2,\al)} \otimes C_{(1,\ali)}.
\end{equation*}
Therefore, for all $\al_1, \dots,\al_n \in \pi$ such that $\al_1 \cdots \al_n=1$, we obtain
\begin{eqnarray*}
\lefteqn{C_{(1,\al_1)} \otimes \cdots \otimes C_{(n-1,\al_{n-1})} \otimes C_{(n,\al_n)}}\\
   & = & (C_{(1,\al_n^{-1})})_{(1,\al_1)} \otimes \cdots \otimes (C_{(1,\al_n^{-1})})_{(n-1,\al_{n-1})}
         \otimes C_{(2,\al_n)} \\
   & = & (C_{(2,\al_n^{-1})})_{(1,\al_1)} \otimes \cdots \otimes (C_{(2,\al_n^{-1})})_{(n-1,\al_{n-1})}
         \otimes C_{(1,\al_n)} \\
   & = & C_{(2,\al_1)} \otimes \cdots \otimes C_{(n,\al_{n-1})} \otimes  C_{(1,\al_n)}.
\end{eqnarray*}
\end{proof}

\section{Invariants of colored Heegaard diagrams}\label{sect-kup1}
In this section, we define \p colored Heegaard diagrams and their equivalence. Then, starting from an involutory Hopf \p
coalgebra, we construct an equivalence invariant of \p colored Heegaard diagrams.

\subsection{Colored Heegaard diagrams}\label{s:colheegdiags}
By a \emph{Heegaard diagram}, we shall mean a triple $D=(S,u,l)$ where $S$ is a closed, connected, and oriented surface
of genus $g\geq 1$ and $u=\{u_1, \dots, u_g \}$ and $l=\{l_1, \dots, l_g \}$ are two systems of pairwise disjoint closed
curves on $S$ such that the complement to $\cup_k u_k$ (resp.\@ $\cup_i l_i$) is connected. Note that if a sphere with
$g$ handles is cut along $g$ disjoint circles that do not split it, then a sphere from which $2g$ disks have been deleted
is obtained.

The circles $u_k$ (resp.\@ $l_i$) are called the \emph{upper} (resp.\@ \emph{lower}) \emph{circles} of the diagram. By
general position we can (and we always do) assume that $u$ and $l$ are transverse. Note that $u \cap l$ is then a finite
set. The Heegaard diagram $D$ is said to be \emph{oriented} if all its lower and upper circles are oriented.

Let $D=(S,u,l)$ be an oriented Heegaard diagram. Denote by $g$ the genus of $S$. Fix an alphabet $X=\{x_1, \dots ,x_g\}$
in $g$ letters. For any $1 \leq i\leq g$, travelling along the lower circle $l_i$ gives a word $w_i(x_1, \dots, x_g)$ as
follows:
\begin{itemize}
 \item Start with the empty word $w_i=\emptyset$;
 \item Make a round trip along $l_i$ following its orientation. Each time $l_i$ encounters an upper circle $u_k$ at some
       crossing $c \in l_i \cap u_k$ (for some $1 \leq k \leq g$), replace $w_i$ by $w_i x_k^\nu$ where:
      \begin{equation*}
        \nu=\begin{cases}
                    +1 & \text{if $(d_c l_i, d_c u_k)$ is an oriented basis for $T_c S$,} \\
                    -1 & \text{otherwise;}
                 \end{cases}
      \end{equation*}
 \item After a complete turn along $l_i$, one gets $w_j$.
\end{itemize}
Note that the word $w_i$ is well-defined up to conjugacy by some word in the letters $x_1, \dots, x_g$ (this is due to
the indeterminacy in the choice of the starting point on $l_i$).

We say that the oriented Heegaard diagram $D$ is \emph{\p colored} if each upper circle $u_k$ is provided with an element
$\al_k \in \pi$ such that $w_i(\al_1, \dots, \al_g)=1 \in \pi$ for all $1 \leq i \leq g$. The system $\al=(\al_1, \dots,
\al_g)$ is called the \emph{color} of $D$.

Two \p colored Heegaard diagrams are said to be \emph{equivalent} if one can be obtained from the other by a finite
sequence of the following moves (or their inverse):

{\bf Type I (homeomorphism of the surface):} By using an orientation-preserving homeomorphism of a (closed, connected,
  and oriented) surface $S$ to a (closed, connected, and oriented) surface $S'$, the oriented upper (resp.\@ lower) circles
  on $S$ are carried to the oriented upper (resp.\@ lower) circles on $S'$. The colors of the upper circles remain unchanged.

{\bf Type II (orientation reversal):} The orientation of an upper or lower circle is changed to its opposite. For an
  upper circle $u_i$, its color $\al_i$ is changed to its inverse $\al_i^{-1}$.

\begin{figure}[h]
                         \psfrag{u}[][]{$u_k$}
                         \psfrag{a}[][]{$\alpha_k$}
                         \psfrag{l}[][]{$l_i$}
                         \scalebox{.9}{\includegraphics{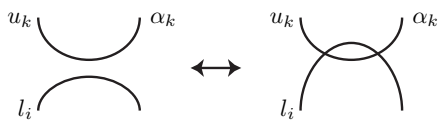}}

     \caption{Two-point move}
     \label{isotopy2}
  \end{figure}

{\bf Type III (isotopy of the diagram):} We isotop the lower circles of the diagram relative to the upper circles. If
  this isotopy is in general position, it reduces to a sequence of two-point moves shown in Figure~\ref{isotopy2}. The
  colors of the upper circles remain unchanged.

\begin{figure}[h]
                        \psfrag{D}[][]{$D$}
                        \scalebox{.9}{\includegraphics{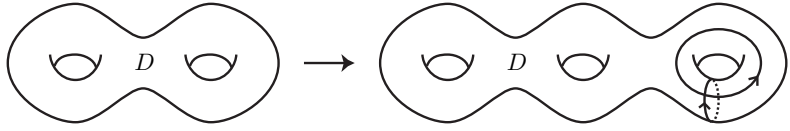}}
     \caption{Stabilization}
     \label{stabilisation2}
  \end{figure}

{\bf Type IV (stabilization):} We remove a disk from $S$ which is disjoint from all upper and lower circles and replace
  it by a punctured torus with one upper and one lower (oriented) circles. One of them corresponds to the standard
  meridian and the other to the standard longitude of the added torus, see Figure~\ref{stabilisation2}. The added upper
  circle is colored with $1 \in \pi$.

\begin{figure}[h]
                         \psfrag{U}[][]{$C'_1$}
                         \psfrag{V}[][]{$C'_2$}
                         \psfrag{b}[][]{$b$}
                         \psfrag{C}[][]{$C_1$}
                         \psfrag{D}[][]{$C_2$}
                         \scalebox{.9}{\includegraphics{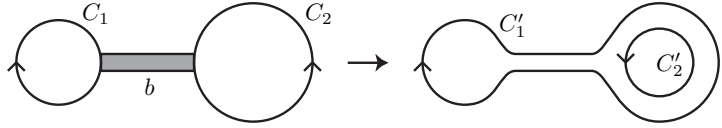}}
     \caption{Circle slide}
     \label{circleslide2}
  \end{figure}

{\bf Type V (sliding a circle past another):} Let $C_1$ and $C_2$ be two circles of a \p colored Heegaard diagram, both
  upper or both lower and let $b$ be a band on $S$ which connects $C_1$ to $C_2$ (that is, $b\co I \times I \to S$ is an embedding
  of $[0,1] \times [0,1]$ for which $b(I \times I)\cap C_i=b(i  \times I)$, $i=1,2$) but does not cross any other
  circle. The circle $C_1$ is replaced by
  \begin{equation*}
  C'_1=C_1 \#_b C_2=C_1 \cup C_2 \cup b(I \times \partial I) \setminus
  b(\partial I \times I).
  \end{equation*}
  The circle $C_2$ is replaced by a copy $C'_2$ of
  itself which is slightly isotoped such that it has no point in common with $C'_1$. The new circle $C'_1$ (resp.\@
  $C'_2$) inherits of the orientation induced by $C_1$ (resp.\@ $C_2$), see Figure~\ref{circleslide2}.

  If the two circles are both lower, then the colors of the upper circle remain unchanged.
  Suppose that the two circles are both upper, say $C_1=u_i$ and $C_2=u_j$ with colors $\al_i$ and $\al_j$ respectively.
  Set $p=b(0,\frac{1}{2}) \in u_i$ and $q=b(1,\frac{1}{2}) \in u_j$.
  Up to first applying a move of type II to $u_i$ and/or $u_j$, we can assume that $(d_p b( \cdot,\frac{1}{2}), d_p u_i)$
  is a negatively-oriented basis for $T_pS$ and $(d_q b( \cdot,\frac{1}{2}), d_q u_j)$ is a positively-oriented basis
  for $T_qS$. Then the color of $u'_i=C'_1$ is $\al_i$
  and the color of $u'_j=C'_2$ is $\al_i^{-1} \al_j$. The colors of the other upper circles remain unchanged.\\

One can remark that all these moves transform a \p colored Heegaard diagram into another \p colored Heegaard diagram.
Indeed, for a move of type I, each word $w_i$ is replaced by a conjugate of itself. For a move of type II applied to an
upper circle $u_k$, each word $w_i(x_1, \cdots, x_k, \cdots, x_g)$ is replaced by a conjugate of $w_i(x_1, \cdots,
x_k^{-1}, \cdots, x_g)$. For a move of type II applied to a lower circle $l_i$, the word $w_i$ is replaced by a conjugate
of $w_i^{-1}$. For a move of type III between $u_k$ and $l_i$, the word $w_i$ is replaced by a conjugate of itself from
which $x_kx_k^{-1}$ or $x_k^{-1}x_k$ has been inserted. For a move of type IV, the new word $w_{g+1}(x_1, \cdots,
x_{g+1})$ is $x_{g+1}^{\pm 1}$. For a move of type V applied to two lower circles, say $l_i$ slides past $l_j$, the word
$w_i$ is replaced by a conjugate of itself from which a conjugate of $w_j^{\pm 1}$ has been inserted and the other words
remain unchanged (up to conjugation). For a move of type V applied to two upper circles, say $u_i$ slides past $u_j$,
each word $w_k(x_1, \cdots, x_j, \cdots, x_g)$ is replaced by a conjugate of $w_k(x_1, \cdots, x_i x_j, \cdots, x_g)$
(see the assumptions on the orientation of the circles $u_i$ and $u_j$). Therefore the conditions $w_i(\al_1, \dots,
\al_g)=1$ are still verified when performing one of these moves.

\subsection{Invariants of \protect\p colored Heegaard diagrams}
Fix a finite type involutory Hopf \p coalgebra $H=\{H_\al\}_{\al \in \pi}$ such that $\dim H_1 \neq 0$ in the ground
field $\kk$ of $H$. Note that $H$ is then semisimple and cosemisimple (by Lemma~\ref{semcosemforinvol}). Using this
algebraic data, we give a method to define an invariant of \p colored Heegaard diagrams, which generalizes that of
Kuperberg~\cite{Ku1}.

Let $D=(S,u,l)$ be a \p colored Heegaard diagram with color $\al=(\al_1,\dots,\al_g)$.

To each upper circle $u_k$, we associate the tensor of Figure~\ref{assignakupderA}, where $c_1, \dots, c_n$ are the
crossings between $u_k$ and $l$ which appear in this order when making a round trip along $u_k$ following its
orientation. Since this tensor is cyclically symmetric (see Lemma~\ref{symabrev}), this assignment does not depend on the
choice of the starting point on the circle $u_k$ when making a round trip along it.
\begin{figure}[h]
       \subfigure[Tensor associated to $u_k$]{\label{assignakupderA}
                         \psfrag{m}[][]{$m_{\al_k}$}
                         \psfrag{c}[][]{$c_1$}
                         \psfrag{a}[][]{$c_n$}
                         \scalebox{.9}{\includegraphics{algokup1der.eps}}}\quad\quad
       \subfigure[Tensor associated to $l_i$]{\label{assignakupderB}
                         \psfrag{u}[][]{$c_m$}
                         \psfrag{c}[][]{$c_1$}
                         \psfrag{d}[][]{$\Delta_{\be_1, \cdots, \be_m}$}
                         \scalebox{.9}{\includegraphics{algokup2der.eps}}}
           \caption{}
           \label{assignakupder}
\end{figure}

To each lower circle $l_i$, we associate the tensor of Figure~\ref{assignakupderB}, where $c_1, \dots, c_m$ are the
crossings between $l_i$ and $u$ which appear in this order when making a round trip along $l_i$ following its
orientation, and the $\be_j\in\pi$ are defined as follows: if $l_i$ intersects $u_k$ at $c_j$, then $\be_j=\al_k^\nu$
with $\nu=+1$ if $(d_{c_j}l_i, d_{c_j}u_k)$ is an oriented basis for $T_{c_j} S$ and $\nu=-1$ otherwise. Note that $\be_1
\cdots \be_m=w_i(\al_1,\dots,\al_g)=1$ and so the tensor associated to $l_i$ is well defined. Since this tensor is
cyclically symmetric (see Lemma~\ref{symabrev}), this assignment does not depend on the choice of the starting point on
the circle $l_i$ when making a round trip along it.

Let $c$ be a crossing point between an upper and a lower circle, say between $u_k$ and $l_i$. Let $\nu$ be as above. If
$\nu=+1$, we contract the tensors assigned to $l_i$ and $u_k$ as follows:
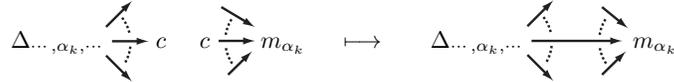
\begin{equation*}
                               \psfrag{d}[B][B]{$\Delta_{\cdots,\al_k,\cdots}$}
                               \psfrag{c}[B][B]{$c$}
                               \psfrag{m}[B][B]{$m_{\al_k}$}
                               \psfrag{=}[B][B]{$\longmapsto$}
                               \scalebox{.9}{\includegraphics{algokupcontrac3.eps}}
\end{equation*}
If $\nu=-1$, we contract the tensors assigned to $l_i$ and $u_k$ as follows:
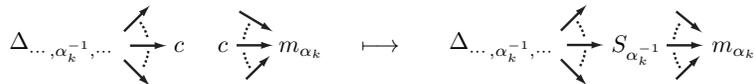
\begin{equation*}
                               \psfrag{d}[B][B]{$\Delta_{\cdots,\al_k^{-1},\cdots}$}
                               \psfrag{c}[B][B]{$c$}
                               \psfrag{=}[B][B]{$\longmapsto$}
                               \psfrag{m}[B][B]{$m_{\al_k}$}
                               \psfrag{S}[B][B]{$S_{\al_k^{-1}}$}
                               \scalebox{.9}{\includegraphics{algokupcontrac4.eps}}
\end{equation*}
After all contractions, one gets $Z(D) \in \Bbbk$.

Finally, we set:
\begin{equation*}
K_H(D)=(\dim H_1 )^{-g} \,Z(D) \in \kk.
\end{equation*}
\begin{theorem}\label{invheegdiaglem1}
Let $H=\{H_\al\}_{\al \in \pi}$ be a finite type involutory Hopf \p coalgebra with $\dim H_1 \neq 0$ in the ground field
$\kk$ of $H$. Then $K_H$ is an equivalence invariant of \p colored Heegaard diagrams.
\end{theorem}
\begin{proof}
We have to verify that $K_H$ is invariant under the moves of type I-V. The proof is similar to that of
\cite[Theorem~5.1]{Ku1}, except that we have to take care of the colors of the Heegaard diagrams.

Clearly, $K_H$ is invariant under a move of type I.

Consider a move of type II applied to an upper $u_k$ circle with color $\al_k$, that is, $u_k$ is replaced by $u'_k=u_k$
with the opposite orientation and with color $\al_k^{-1}$. Let $c_1, \dots, c_n$ be the crossings between $u_k$ and the
lower circles which appear in this order following the orientation. Then the tensor associated to $u_k$ (resp.\@ $u'_k$)
is:
\begin{equation*}
                               \psfrag{m}[][]{$m_{\al_k}$}
                               \psfrag{c}[][]{$c_1$}
                               \psfrag{n}[][]{$m_{\al_k^{-1}}$}
                               \psfrag{a}[][]{$c_n$}
                               \psfrag{r}[][]{$\Bigl($ resp.}
                               \psfrag{p}[][]{$\Bigr ).$}
                               \scalebox{.9}{\includegraphics{invarkup1.eps}}
\end{equation*}
Recall that the contraction rule applied to a crossing point $c \in u_k \cap l_i$ is:
\begin{equation*}
                               \psfrag{d}[][]{$\Delta_{\cdots,\al_k^\nu,\cdots}$}
                               \psfrag{p}[][]{$\phi$}
                               \psfrag{a}[][]{$c_1$}
                               \psfrag{m}[][]{$m_{\al_k}$}
                               \psfrag{b}[][]{$c_n$}
                               \scalebox{.9}{\includegraphics{invarkup2.eps}}
\end{equation*}
where $\nu=+1$ and $\phi=\id_{H_{\al_k}}$ if $(d_cl_i,d_c u_k)$ is a positively-oriented basis of $T_cS$ and $\nu=-1$ and
$\phi=S_{\al_k^{-1}}$ otherwise. Then the contraction rule applied to the corresponding crossing point $c' \in u'_k \cap
l_i$ is:
\begin{equation*}
                               \psfrag{d}[][]{$\Delta_{\cdots,\al_k^\nu,\cdots}$}
                               \psfrag{p}[][]{$\psi$}
                               \psfrag{a}[][]{$c_n$}
                               \psfrag{m}[][]{$m_{\al_k^{-1}}$}
                               \psfrag{b}[][]{$c_1$}
                               \scalebox{.9}{\includegraphics{invarkup2b.eps}}
\end{equation*}
where $\psi=\id_{H_{\al_k^{-1}}}$ if $(d_cl_i,d_c u'_k)$ is a positively-oriented basis of $T_cS$ and $\psi=S_{\al_k}$
otherwise. Now $\psi= S_{\al_k}\circ \phi$ since the antipode is involutory. Therefore the invariance follows from the
equality:
\begin{equation*}
                               \psfrag{=}[][]{$=$}
                               \psfrag{S}[][]{$S_{\al_k}$}
                               \psfrag{c}[][]{$c_1$}
                               \psfrag{m}[Bc][Bc]{$m_{\al_k}$}
                               \psfrag{n}[Bc][Bc]{$m_{\al_k^{-1}}$}
                               \psfrag{a}[][]{$c_n$}
                               \scalebox{.9}{\includegraphics{invarkup3.eps}}
\end{equation*}
which comes from the anti-multiplicativity of the antipode (see \cite[Lemma~1.1(a)]{Vir2}) and the fact that $T_\ali
\circ S_\al = T_\al$ for any $\al \in \pi$ (by Lemma~\ref{TCtwosided}).

For a move of type II applied to a lower circle, the invariance follows from the equality:
\begin{equation*}
                               \psfrag{=}[][]{$=$}
                               \psfrag{S}[][]{$S_{\be_1^{-1}}$}
                               \psfrag{T}[][]{$S_{\be_m^{-1}}$}
                               \psfrag{u}[][]{$c_1$}
                               \psfrag{d}[][]{$\Delta_{\be_1,\cdots,\be_m}$}
                               \psfrag{f}[][]{$\Delta_{\be_m^{-1},\cdots,\be_1^{-1}}$}
                               \psfrag{c}[][]{$c_m$}
                               \scalebox{.9}{\includegraphics{invarkup4.eps}}
\end{equation*}
which comes from the anti-comultiplicativity of the antipode (see \cite[Lemma~1.1(c)]{Vir2}) and the fact that $S_1(C) =
C$ (by Lemma~\ref{TCtwosided}).

Consider now a two-point move between an upper circle with color $\al$ and a lower circle. Up to first applying a move of
type II, we can consider that these two circles are oriented so that the invariance is a consequence of the following
equality:
\begin{equation*}
                               \psfrag{=}[][]{$=$}
                               \psfrag{S}[][]{$S_\ali$}
                               \psfrag{m}[][]{$m_\al$}
                               \psfrag{d}[][]{$\Delta_{\cdots,\al,\ali, \cdots}$}
                               \psfrag{f}[][]{$\Delta_{\cdots,1, \cdots}$}
                               \psfrag{g}[][]{$\Delta_{\cdots, \cdots}$}
                               \psfrag{e}[][]{$\varepsilon$}
                               \psfrag{i}[][]{$1_\al$}
                               \scalebox{.9}{\includegraphics{invarkup5.eps}}
\end{equation*}
which comes from \eqref{antipode}.

A move of type IV contributes $C \to T_1=\dim H_1$ (see Lemma~\ref{TCtwosided}) to $Z(D)$, which is cancelled by
normalization.

Consider a move of type V applied to two upper circles, say $u_i$ (with color $\al_i$) slides past $u_j$ (with color
$\al_j$). We assume, as a representative case, that both circles have three crossings with the lower circles:
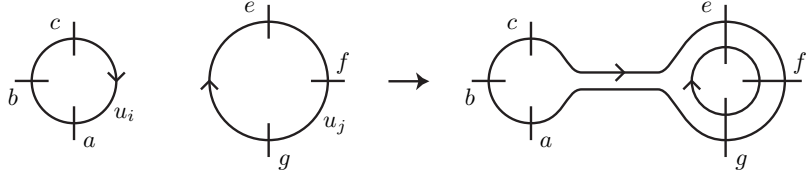
\begin{equation*}
                        \psfrag{a}[][]{$a$}
                        \psfrag{b}[][]{$b$}
                        \psfrag{c}[][]{$c$}
                        \psfrag{e}[][]{$e$}
                        \psfrag{f}[][]{$f$}
                        \psfrag{g}[][]{$g$}
                        \psfrag{u}[][]{$u_i$}
                        \psfrag{v}[][]{$u_j$}
                        \scalebox{.9}{\includegraphics{invarkup6.eps}}
\end{equation*}
Using the anti-multiplicativity of the antipode (which allows us to consider only the positively-oriented case of the
contraction rule), we have that the following factor of $Z(D)$:
\begin{equation*}
                        \psfrag{a}[][]{$a$}
                        \psfrag{b}[][]{$b$}
                        \psfrag{c}[][]{$c$}
                        \psfrag{e}[][]{$e$}
                        \psfrag{f}[][]{$f$}
                        \psfrag{g}[][]{$g$}
                        \psfrag{n}[Bc][Bc]{$m_{\al_j}$}
                        \psfrag{m}[Bc][Bc]{$m_{\al_i}$}
                        \scalebox{.9}{\includegraphics{invarkup7.eps}}
\end{equation*}
is replaced by:
\begin{equation*}
                        \psfrag{a}[][]{$a$}
                        \psfrag{b}[][]{$b$}
                        \psfrag{c}[][]{$c$}
                        \psfrag{e}[][]{$e$}
                        \psfrag{f}[][]{$f$}
                        \psfrag{g}[][]{$g$}
                        \psfrag{n}[][]{$m_{\al_i^{-1}\al_j}$}
                        \psfrag{m}[][]{$m_{\al_i}$}
                        \psfrag{D}[][]{$\cp{\al_i}{\al_i^{-1}\al_j}$}
                        \scalebox{.9}{\includegraphics{invarkup8.eps}}
\end{equation*}
By using the multiplicativity of the comultiplication and the fact that $(T_\al)_{\al \in \pi}$ is a left \p integral for
$H$ (see Lemma~\ref{TCtwosided}), we obtain that these two factors are equal, as depicted in Figure~\ref{finavnt13dec1}.
\begin{figure}[h]
                        \psfrag{=}[][]{$=$}
                        \psfrag{a}[][]{$a$}
                        \psfrag{b}[][]{$b$}
                        \psfrag{c}[][]{$c$}
                        \psfrag{e}[][]{$e$}
                        \psfrag{f}[][]{$f$}
                        \psfrag{g}[][]{$g$}
                        \psfrag{n}[][]{$m_{\al_i^{-1}\al_j}$}
                        \psfrag{m}[][]{$m_{\al_i}$}
                        \psfrag{w}[][]{$m_{\al_j}$}
                        \psfrag{i}[][]{$1_{\al_i}$}
                        \psfrag{R}[][]{$T_{\al_i^{-1}\al_j}$}
                        \psfrag{T}[][]{$T_{\al_i}$}
                        \psfrag{L}[][]{$T_{\al_j}$}
                        \psfrag{D}[][]{$\cp{\al_i}{\al_i^{-1}\al_j}$}
                        \scalebox{.9}{\includegraphics{invarkup9.eps}}
     \caption{}
     \label{finavnt13dec1}
\end{figure}

Finally, suppose that a lower circle slides past another lower circle. We assume, as a representative case, that these
two circles have both three crossings with the upper circles. Let $\al_1,\al_2,\al_3$ (resp.\@ $\be_1,\be_2,\be_3$) be
the colors of the upper circles intersected (following the orientation) by the first (resp. second) lower circle
considered. Then the invariance follows from the equality of Figure~\ref{finavnt13dec2}
\begin{figure}[t]
                        \psfrag{=}[][]{$=$}
                        \psfrag{x}[][]{$\Delta_{1,1}$}
                        \psfrag{n}[][]{$m_{\be_2}$}
                        \psfrag{m}[][]{$m_{\be_1}$}
                        \psfrag{w}[][]{$m_{\be_3}$}
                        \psfrag{e}[][]{$\varepsilon$}
                        \psfrag{u}[][]{$m_1$}
                        \psfrag{C}[][]{$C$}
                        \psfrag{q}[][]{$\Delta_{\be_1,\be_2,\be_3}$}
                        \psfrag{r}[][]{$\Delta_{\al_1,\al_2,\al_3}$}
                        \psfrag{d}[][]{$\Delta_{\al_1,\al_2,\al_3,\be_1,\be_2,\be_3}$}
                        \scalebox{.9}{\includegraphics{invarkup10.eps}}
     \caption{}
     \label{finavnt13dec2}
\end{figure}
which comes from the multiplicativity of the comultiplication and the fact that $C$ is a right integral for $H_1$ (see
Lemma~\ref{TCtwosided}). This completes the proof of the theorem.
\end{proof}

If $D$ is a \p colored Heegaard diagram of genus $g$ with color $\al=(\al_1, \dots, \al_g)$ and $\be \in \pi$, then
$\be\al\bei=(\be\al_1\bei, \dots, \be\al_g\bei)$ is clearly another color of the Heegaard diagram, said to be
\emph{conjugate} to the color $\alpha$.

\begin{lemma}\label{cojugclassindepKH}
Let $H=\{H_\al\}_{\al \in \pi}$ be a finite type involutory Hopf \p coalgebra with $\dim H_1 \neq 0$ in the ground field
$\kk$ of $H$. Suppose that $H=\{H_\al\}_{\al \in \pi}$ is crossed (see Section~\ref{deficro}). Then the invariant
$K_H(D)$ does not depend on the conjugacy class of the color of the \p colored Heegaard diagram~$D$.
\end{lemma}
\begin{proof}
Suppose that $H=\{H_\al\}_{\al \in \pi}$ admits a crossing $\varphi=\{\varphi_\be\co H_\al \to H_{\be \al \bei}
\}_{\al,\be \in \pi}$. Let $D=(S,u,l)$ be a \p colored Heegaard diagram of genus $g$ with color $\al=(\al_1, \dots,
\al_g)$. Fix $\be \in \pi$. Let us denote by $D^\be$ the underlying oriented Heegaard diagram of $D$ endowed with the
color $\be\al\bei=(\be\al_1\bei, \dots, \be\al_g\bei)$. We have to verify that $K_H(D^\be)=K_H(D)$.

Let $1 \leq k,i \leq g$ and denote by $c_1, \dots, c_n$ (resp.\@ $c'_1, \dots, c'_m$) the crossings between $u_k$ and $l$
(resp.\@ $l_i$ and $u$) which appear in this order when making a round trip along $u_k$ (resp.\@ $l_i$) following its
orientation. Recall that, for $D$ (resp.\@ $D^\be$), the tensor of Figure~\ref{lemheegconj1a} (resp.\@
Figure~\ref{lemheegconj1b}) is associated to the upper circle $u_k$, and the tensor of Figure~\ref{lemheegconj1c}
(resp.\@ Figure~\ref{lemheegconj1d}) is associated to the lower circle $l_i$ where, if $l_i$ intersects some $u_n$ at
$d_j$, $\be_j=\al_n^\nu$ with $\nu=1$ if $(d_{d_j}l_i, d_{d_j}u_n)$ is an oriented basis for $T_{d_j} S$ and $\nu=-1$
otherwise.
\begin{figure}[h]
   \subfigure[]{\label{lemheegconj1a}
                         \psfrag{m}{$m_{\al_k}$}
                         \psfrag{c}[Br][Br]{$c_1$}
                         \psfrag{a}[Br][Br]{$c_n$}
                         \scalebox{.9}{\includegraphics{conjkup1.eps}}} \hspace*{.5cm}
   \subfigure[]{\label{lemheegconj1b}
                         \psfrag{m}{$m_{\be\al_k\bei}$}
                         \psfrag{c}[Br][Br]{$c_1$}
                         \psfrag{a}[Br][Br]{$c_n$}
                         \scalebox{.9}{\includegraphics{conjkup1b.eps}}} \hspace*{.8cm}
   \subfigure[]{\label{lemheegconj1c}
                         \psfrag{u}{$c'_m$}
                         \psfrag{c}{$c'_1$}
                         \psfrag{d}[Br][Br]{$\Delta_{\be_1, \cdots, \be_m}$}
                         \scalebox{.9}{\includegraphics{conjkup2.eps}}} \hspace*{.8cm}
   \subfigure[]{\label{lemheegconj1d}
                         \psfrag{u}{$c'_m$}
                         \psfrag{c}{$c'_1$}
                         \psfrag{d}[Br][Br]{$\Delta_{\be\be_1\bei, \cdots, \be\be_m\bei}$}
                         \scalebox{.9}{\includegraphics{conjkup2b.eps}}}
     \caption{}
     \label{lemheegconj01}
\end{figure}

Since $H$ is cosemisimple, \cite[Lemmas~6.3(a) and 7.2]{Vir2} give that $\varphi_\be(C)=C$ and \cite[Corollary~6.2  and
Lemma 7.2]{Vir2} give that $T_{\be \al \bei} \pb= T_\al$ for all $\al\in \pi$. Therefore, using the multiplicativity and
comultiplicativity of a crossing, we have the equalities of Figure~\ref{lemheegconj2}.
\begin{figure}[h]
                         \psfrag{p}{$\pbi$}
                         \psfrag{q}{$\pb$}
                         \psfrag{n}{$m_{\al_k}$}
                         \psfrag{m}{$m_{\be\al_k\bei}$}
                         \psfrag{u}[Br][Br]{$\Delta_{\be_1, \cdots, \be_m}$}
                         \psfrag{d}[Br][Br]{$\Delta_{\be\be_1\bei, \cdots, \be\be_m\bei}$}
                         \psfrag{=}{$=$}
                         \scalebox{.9}{\includegraphics{conjkup3.eps}}
     \caption{}
     \label{lemheegconj2}
\end{figure}

Hence, since $\varphi_\be \circ \varphi_\bei=\id_{H_\al}$ and $S_{\be \al \bei} \varphi_\be = \varphi_\be S_\al$ for all
$\al \in \pi$ by \cite[Lemma~6.1]{Vir2}, contracting the tensors  associated to $D^\be$ and $D$ leads to the same scalar
$Z(D^\be)=Z(D)$. Finally $K_H(D^\be)=(\dim H_1 )^{-g} Z(D^\be)=(\dim H_1 )^{-g} Z(D)=K_H(D)$.
\end{proof}

\section{Invariants of flat bundles over 3-manifolds}\label{sect-kup2}

In this section, we use the invariant of colored Heegaard diagrams defined in Section~\ref{sect-kup1} to construct an
invariant of flat bundles over \trois manifolds.

\subsection{\protect Flat bundles over 3-manifolds}\label{triplet}
Fix a group $\pi$. By a \emph{flat \p bundle over a \trois manifold}, we shall mean a principal \p bundle $\xi=(p\co
\Tilde{M} \to M)$, where $M$ is a closed connected and oriented \trois manifold, which is flat, that is, such that its
transition functions are locally constant. Such an object can be viewed as a regular covering $\Tilde{M} \to M$ with
group of automorphisms $\pi$. The space $\Tilde{M}$ (resp.\@ $M$) is called the \emph{total space} (resp.\@ \emph{base
space}) of $\xi$.

Two flat \p bundles over \trois manifolds $\xi$ and $\xi'$ are said to be \emph{equivalent} if there exists an
homeomorphism $\tilde{h}\co \Tilde{M} \to \Tilde{M'}$ between their total spaces which preserves the action of $\pi$ and
which induces an orientation-preserving homeomorphism $h\co M \to M'$ between their base spaces.

A flat \p bundle $\xi=(p\co \Tilde{M} \to M)$ is said to be \emph{pointed} when its total space $\Tilde{M}$ is endowed
with a base point $\tilde{x} \in \Tilde{M}$. Two pointed flat \p bundles over \trois manifolds $(\xi,\tilde{x})$ and
$(\xi',\tilde{x}')$ are said to be \emph{equivalent} if there exists an equivalence $\tilde{h}\co \Tilde{M} \to
\Tilde{M'}$ between them such that $\tilde{h}(\tilde{x})=\tilde{x}'$.

Let $(\xi,\tilde{x})$ be a pointed flat \p bundle over a \trois manifold. Set $x=p(\tilde{x})\in M$, where $p\co
\Tilde{M} \to M$ is the bundle map of $\xi$. We can associate to $(\xi,\tilde{x})$ a morphism $f\co \pi_1(M,x) \to \pi$,
called \emph{monodromy of $\xi$ at $\tilde{x}$}, by the following procedure: any loop $\ga$ in $(M,x)$ uniquely lifts to
a path $\tilde{\ga}$ in $\Tilde{M}$ beginning at $\tilde{x}$. The path $\tilde{\ga}$ ends at $\al \cdot \tilde{x}$ for a
unique $\al \in \pi$. The monodromy is defined by $f([\ga])=\al$, where $[\ga]$ denotes the homotopy class in
$\pi_1(M,x)$ of the loop~$\ga$.

Any pointed flat \p bundle $(\xi,\tilde{x})$ over a \trois manifold $M$ leads to a triple $(M,x,f)$, where $x$ is the
image of $\tilde{x}$ under the bundle map $\Tilde{M} \to M$ of $\xi$ and $f\co \pi_1(M,x) \to \pi$ is the monodromy of
$\xi$ at $\tilde{x}$. Conversely, a triple $(M,x,f)$, where $M$ is a (closed, connected,  and oriented) \trois manifold,
$x \in M$, and $f\co \pi_1(M,x) \to \pi$ is a group homomorphism, leads to a pointed flat \p bundle over $M$ uniquely
determined up to equivalence (see \cite[Proposition~14.1]{Ful}). Let us briefly recall the construction of the pointed
flat \p bundle $\xi=(p\co\Tilde{M} \to M)$ associated to a triple $(M,x,f)$: the pointed manifold $(M,x)$ admits a
universal covering $u\co (Y,y) \to (M,x)$. The fundamental group $\pi_1(M,x)$ acts on the left on $Y$ as follows: given
$[\sigma]\in \pi_1(M,x)$ and $z\in Y$, let $y'$ be the endpoint of the lift (in $Y$) of $\sigma$ to a path that starts on
$y$. Choose a path $\gamma$ from $y$ to $z$ in $Y$. Then $[\sigma]\cdot z$ is defined to be the endpoint of the lift (in
$Y$) of the path $u\circ \gamma$ that starts at $y'$. Set $\Tilde{M}=(Y \times \pi) /\pi_1(M,x)$, where $\pi_1(M,x)$ acts
on the left on $Y \times \pi$ by the rule $[\sigma] \cdot (z,\alpha)=([\sigma]\cdot z, \alpha f([\sigma])^{-1})$. Let
$\langle z,\alpha \rangle \in \Tilde{M}$ denote the image of $(z,\alpha) \in Y \times \pi$ and define $p\co  \Tilde{M}
\to M$ by $p(\langle z,\alpha \rangle)=u(z)$. Then $\xi=(p\co \Tilde{M} \to M)$ is a flat \p bundle over $M$ which is
pointed with base point $\tilde{x}=\langle y,1 \rangle$ and whose monodromy at $\tilde{x}$ is $f$.

It may be convenient to adopt this second point of view. In particular, under this point of view, two pointed flat \p
bundles over \trois manifolds $(M,x,f)$ and $(M',x',f')$ are equivalent if there exists an orientation-preserving
homeomorphism $h\co M \to M'$ such that $h(x)=x'$ and $f' \circ h_* = f$, where $h_*\co \pi_1(M,x) \to \pi_1(M',x')$ is
the group isomorphism induced by $h$ in homotopy.\\

\subsection{Heegaard diagrams of 3-manifolds}
We first recall that a \emph{Heegaard splitting of genus $g$} of a closed, connected, and oriented \trois manifold $M$ is
an ordered pair $(M_u,M_l)$ of submanifolds of $M$, both homeomorphic to a handelbody of genus $g$, such that $M=M_u \cup
M_l$ and $M_u\cap M_l=\partial M_u=\partial M_l$. The handelbody $M_u$ (resp.\@ $M_l$) is called \emph{upper} (resp.\@
\emph{lower}) and the surface $\partial M_u=\partial M_l$ is called a \emph{Heegaard surface} (of genus $g$) for $M$.

It is well known that every closed, connected, and oriented \trois manifold $M$ has a Heegaard splitting (e.g., by taking
a closed regular neighborhood of the one-dimensional skeleton of a triangulation of $M$ and the closure of its
complement).

Let $(M_u,M_l)$ be a Heegaard splitting of genus $g$ of a closed, connected, and oriented \trois manifold $M$. Since
$M_u$ is homeomorphic to a handelbody of genus $g$, there exists a finite collection $\{D_1, \cdots , D_g\}$ of pairwise
disjoint properly embedded $2$-disks in $M_u$ which cut $M_u$ into a $3$-ball. Likewise, there exists a finite collection
$\{D'_1, \cdots , D'_g\}$ of pairwise disjoint properly embedded $2$-disks in $M_l$ which cut $M_l$ into a $3$-ball. For
$1 \leq i \leq g$, set $u_i=\partial D_i$ and $l_i=\partial D'_i$. We can (and we do) suppose that these circles meet
transversely. Denote the Heegaard surface $M_u\cap M_l$ by $S$. It is oriented as follows: for any point $p\in S$, a
basis $(e_1,e_2)$ of $T_p S$ is positive if, when completing $(e_1,e_2)$ with a vector $e_3$ pointing from $M_l$ to
$M_u$, we obtain a positively-oriented a basis $(e_1,e_2,e_3)$ of $T_pM$. Then $D=(S,u=\{u_1, \cdots , u_g\},l=\{l_1,
\cdots , l_g\})$ is a Heegaard diagram in the sense of Section~\ref{s:colheegdiags}. Such a Heegaard diagram is called a
\emph{Heegaard diagram (of genus $g$) of $M$}.

\subsection{Kuperberg-type invariants of flat bundles over 3-manifolds}
Fix a finite type involutory Hopf \p coalgebra $H=\{H_\al\}_{\al \in \pi}$ with $\dim H_1\neq 0$ in the ground field
$\kk$ of $H$.

Let $(\xi=(p\co \Tilde{M} \to M),\tilde{x})$ be a poin\-ted flat \p bundle over a 3-manifold $M$. Set $x=p(\tilde{x})\in
M$ and let $f\co \pi_1(M,x) \to \pi$ be the monodromy of $\xi$ at $\tilde{x}$. Consider a Heegaard diagram $D=(S,u,l)$ of
genus $g$ of $M$. Recall that $S=\partial M_u=\partial M_l \subset M$ for some Heegaard splitting $(M_u,M_l)$ of $M$. We
arbitrarily orient the upper and lower circles so that $D$ is oriented. We can (and we always do) assume that $x \in S
\setminus \{u,l\}$.

Since $S \setminus u$ is homeomorphic to a sphere from which $2g$ disks have been deleted, there exists $g$ pairwise
disjoint (except in $x$) loops $\ga_1, \dots, \ga_g$ on $(S,x)$ such that, for any $1 \leq i \leq g$,
\begin{itemize}
  \item $\ga_i$ intersects the upper circle $u_i$ in exactly one point $p_i$ in such a way that
        $(d_{p_i} \ga_i, d_{p_i} u_i)$ is a positively-oriented basis of $T_{p_i}S$;
  \item $\ga_i$ does not intersect any other upper circle.
\end{itemize}
Then the homotopy classes $a_i=[\ga_i] \in \pi_1(M,x)$ do not depend on the choice of the loops $\ga_i$ verifying the
above conditions (since each $\ga_i$ is homotopic to a unique leaf of the $x$-based $g$-leafed rose formed by the core of
the handelbody $M_u$).

\begin{lemma}\label{VKT}
$\pi_1(M,x)=\langle a_1, \dots, a_g \, | \, w_1(a_1, \cdots, a_g)=1, \dots, w_g(a_1, \cdots, a_g)=1 \rangle$, where the
words $w_i(x_1, \cdots, x_g)$ are defined as in Section~\ref{s:colheegdiags}.
\end{lemma}
\begin{proof}
Recall that there exists a finite collection $\{D_1, \cdots , D_g\}$ (resp.\@ $\{D'_1, \cdots , D'_g\}$) of pairwise
disjoint properly embedded $2$-disks in $M_u$ (resp.\@ $M_l$) which cut $M_u$ (resp.\@~$M_l$) into a $3$-ball and such
that $u_i=\partial D_i$ and $l_i=\partial D'_i$ for $1 \leq i \leq g$. Since $\cup_{1 \leq i \leq g} D_i\cup S$ is a
deformation retract of $M_u$ from which some \trois balls have been deleted, $\pi_1(\cup_{1 \leq i \leq g} D_i\cup S,x)$
is the free group generated by the homotopy classes of the loops $\ga_1, \dots, \ga_g$. Now, by the Seifert-Van Kampen
Theorem, gluing a disk $D'_j$ amounts adding the relation $w_j([\ga_1], \dots,[\ga_g])=1$. Hence the lemma follows from
the fact that $\cup_{1 \leq i \leq g} (D_i \cup D'_i) \cup S$ is a deformation retract of $M$ from which some \trois
balls have been deleted.
\end{proof}

For any $1 \leq i \leq g$, set $\al_i=f(a_i) \in \pi$. By Lemma~\ref{VKT}, $\al=(\al_1, \cdots, \al_g)$ is a color of the
oriented Heegaard diagram $D$. We say that the (oriented) Heegaard diagram $D$ of $M$ is \emph{colored by $f$}.

Finally, we set:
\begin{equation*}
K_H(\xi,\tilde{x})=K_H(D)\in \kk,
\end{equation*}
where $K_H$ is the invariant of \p colored Heegaard diagrams of Theorem~\ref{invheegdiaglem1}.

\begin{theorem}\label{invKuth1}
Let $H=\{H_\al\}_{\al \in \pi}$ be a finite type involutory Hopf \p coalgebra with $\dim
H_1\neq 0$ in the ground field $\kk$ of $H$.
\begin{enumerate}
\renewcommand{\labelenumi}{{\rm (\alph{enumi})}}
\item $K_H$ is an invariant of pointed flat \p bundles over \trois manifolds.
\item Let $\xi=(p\co \Tilde{M} \to M)$ be a flat \p bundle over a \trois manifold $M$.
      \begin{itembul}{{\rm (ii)}}
\item[{\rm (i)}] The function $\tilde{x} \in \Tilde{M}
\mapsto K_H(\xi,\tilde{x}) \in \kk$ is constant on the path-connected components of~$\Tilde{M}$;
\item[{\rm (ii)}] If $H$ is crossed or if $\pi$ is abelian or if the monodromy of $\xi$ is surjective,
then $K_H(\xi,\tilde{x})$ does not depend on the choice
of the base point $\tilde{x} \in \Tilde{M}$.
\end{itembul}
\end{enumerate}
\end{theorem}

This theorem is proved in Section~\ref{preuvethm}.

If $\pi=1$ and $M$ is a \trois manifold, then $K_H(\id_M\co  M \to M)$ coincides with the invariant of $M$ constructed by
Kuperberg \cite{Ku1}.

In Section~\ref{sect-kup3}, we give examples which show that the invariant $K_H$ is not trivial.

\subsection{Basic properties of $K_H$}
Let $H=\{H_\al\}_{\al \in \pi}$ be a finite type involutory Hopf \p coalgebra with $\dim H_1\neq 0$ in the ground field
$\kk$ of $H$. Recall that $H^\opp$ and  $H^\cop$ denotes the opposite or coopposite Hopf \p coalgebra to $H$ (see
Section~\ref{hopfpicoal}).

Let $(\xi,\tilde{x})$ be a pointed flat \p bundle over a \trois manifold $M$. Denote by $-\xi$ the flat \p bundle $\xi$
whose base space $M$ is endowed with the opposite orientation. Then
\begin{equation}\label{eqopmanheeg}
K_H(-\xi,\tilde{x}) =K_{H^\cop}(\xi,\tilde{x})=K_{H^\opp}(\xi,\tilde{x}).
\end{equation}
Indeed, starting from an oriented Heegaard diagram $D=(S,u,l)$ for $M$, reversing the orientation of $M$ amounts
reversing the orientation of the Heegaard surface $S$, and so the first equality of~\eqref{eqopmanheeg} can be easily
obtained by reversing the orientation of the lower circles and the second one by reversing the orientation of the upper
circles.

Let $(\xi_1,\tilde{x}_1)$ and $(\xi_2,\tilde{x}_2)$ be two pointed flat \p bundles over \trois manifolds. Denote by $x_1$
(resp.\@ $x_2$) the image of $\tilde{x}_1$ (resp.\@ $\tilde{x}_2$) under the bundle map of $\Tilde{M}_1 \to M_1$ of
$\xi_1$ (resp.\@ $\Tilde{M}_2 \to M_2$ of $\xi_2$), and by $f_1\co \pi_1(M_1) \to \pi$ (resp.\@ $f_2\co \pi_1(M_2) \to
\pi$) the monodromy of $\xi_1$ at $\tilde{x}_1$ (resp.\@ of $\xi_2$ at $\tilde{x}_2$). Take closed \trois balls $B_1
\subset M_1$ and $B_2 \subset M_2$ such that $x_1 \in \partial B_1$ and $x_2 \in\partial B_2$. Glue $M_1 \setminus {\rm
Int}B_1$ and $M_2 \setminus {\rm Int}B_2$ along a homeomorphism $h\co  \partial B_1 \to \partial B_2$ chosen so that
$h(x_1)=x_2$ and that the orientations in $M_1 \setminus {\rm Int}B_1$ and $M_2 \setminus {\rm Int}B_2$ induced by those
in $M_1$, $M_2$ are compatible. This gluing yields a closed, connected, and oriented \trois manifold $M_1 \# M_2$. For
$i=1$ or $2$, consider the embeddings  $j_i\co  M_i \setminus {\rm Int}B_i \hookrightarrow M_i$ and $k_i\co M_i \setminus
{\rm Int}B_i \hookrightarrow M_1 \# M_2$ and set $x=k_1(x_1)=k_2(x_2)$. By the Van Kampen theorem, since $ \partial B_2
\cong h(\partial B_1)$ is simply-connected, there exists a unique group homomorphism $f\co \pi_1(M_1 \# M_2,x) \to \pi$
such that $f \circ (k_i)_*=f_i \circ (j_i)_*$ ($i=1,2$). This leads to a triple $(M_1 \# M_2,x,f)$, and so to a pointed
flat \p bundle over $M_1 \# M_2$, denoted by $(\xi_1 \# \xi_2, \tilde{x})$, whose monodromy is~$f$ (see
Section~\ref{triplet}). Then
\begin{equation}
K_H(\xi_1 \# \xi_2, \tilde{x})=K_H(\xi_1,\tilde{x}_1) \, K_H(\xi_2,\tilde{x}_2).
\end{equation}
Indeed we can choose a Heegaard diagram for $M$ which is a connected sum of Heegaard diagrams for $M_1$ and $M_2$ and
such that the colorations of these diagrams by the monodromies $f$, $f_1$, or $f_2$ are compatible with this connected
sum.

\subsection{Proof of Theorem~\ref{invKuth1}}\label{preuvethm}

Let us prove Part (a) of Theorem~\ref{invKuth1}. Adopting the second point of view of Section~\ref{triplet}, let
$(M,x,f)$ and $(M',x',f')$ be two equivalent pointed flat \p bundles over \trois manifolds. Let $D$ (resp.\@ $D'$) be an
oriented Heegaard diagrams of $M$ (resp.\@ $M'$) colored by $f$ (resp.\@ $f'$). By virtue of
Theorem~\ref{invheegdiaglem1}, it suffices to prove that $D$ and $D'$ are equivalent \p colored Heegaard diagrams, i.e.,
that $D$ can be obtained from $D'$ by a finite sequence of the moves of type I-V (or their inverses) described in
Section~\ref{s:colheegdiags}.

Since $(M,x,f)$ and $(M',x',f')$ are equivalent, there exists an orientation-pre\-ser\-ving homeomorphism $h \co  M \to
M'$ with $f(x)=x'$ and $f=f' \circ h_*$, where $h_*\co \pi_1(M,x) \to \pi_1(M',x')$ is the homomorphism induced by $h$ in
homotopy. By the Reidemeister-Singer Theorem (see \cite[Theorem~8]{Si} or \cite[Theorem~4.1]{Ku1}), there exist:
\begin{itembul}{\quad $\bullet$}
  \item[\quad $\bullet$] a finite sequence $M_0=M$, $M_1$, $\dots$, $M_{n-1}$, $M_n=M'$ of closed, connected, and oriented
                   \trois manifolds;
  \item[\quad $\bullet$] a Heegaard diagram $D_k=(S_k,u^k=\{u^k_1, \cdots, u^k_{g_k}\},l^k=\{l^k_1, \cdots, l^k_{g_k}\})$
                   of genus $g_k$ of $M_k$ for each $0 \leq k \leq n$, with $D_0=D$ and $D_n=D'$;
  \item[\quad $\bullet$] a finite sequence $h_1\co M_0 \to M_1$, $\dots$,
                   $h_n\co M_{n-1} \to M_n$ of orientation-preserving homeomorphisms;
\end{itembul}
such that $h=h_n \circ \cdots \circ h_1$ and, for any $1 \leq k \leq n$, the diagrams $D_{k-1}$ and $D_k$ are related by
a move (or its inverse) of the following type: \vspace*{.3em}
\begin{itembul}{\quad {\bf Type A:}}
  \item[\quad {\bf Type A:}] $S_k=h_k(S_{k-1})$, $u^k=h_k(u^{k-1})$, and $l^k=h_k(l^{k-1})$;\vspace*{.3em}

\item[\quad {\bf Type B:}] $S_k=h_k(S_{k-1})$, $u^k=h_k(u^{k-1})$, and $l^k$ is isotopic to $h_k(l^{k-1})$ relative to
  $u^k$;\vspace*{.3em}

\item[\quad {\bf Type C:}] $S_k=h_k(S_{k-1}) \# T^2$, $u^k=h_k(u^{k-1}) \cup \{C_1\}$, and
                $l^k=h_k(l^{k-1}) \cup\{C_2\}$, where
                $T^2 $ is a torus and $\{C_1,C_2\}$ is the set formed by the standard meridian and
                longitude of $T^2$;\vspace*{.3em}

\item[\quad {\bf Type D:}] $S_k=h_k(S_{k-1})$, $u^k=h_k(u^{k-1})$, and $l^k$ is obtained from $h_k(l^{k-1})$
           by sliding one circle of
           $h_k(l^{k-1})$ past another circle of $h_k(l^{k-1})$, avoiding the other upper and lower circles of
           $h_k(S_{k-1})$;\vspace*{.3em}

\item[\quad {\bf Type E:}] $S_k=h_k(S_{k-1})$, $l^k=h_k(l^{k-1})$, and $u^k$ is obtained from $h_k(u^{k-1})$
           by sliding one circle of
           $h_k(u^{k-1})$ past another circle of $h_k(u^{k-1})$, avoiding the other upper and lower circles of
           $h_k(S_{k-1})$.\vspace*{.3em}
\end{itembul}

Set $x_0=x \in M_0$ and define $x_k=h_k \circ \cdots h_1(x) \in M_k$ for any $1 \leq k \leq n$. Note that $x_n=x'$ since
$h(x)=x'$. Without loss of generality, we can assume that $x_k \in S_k \setminus \{u^k,l^k\}$. Set $f_0=f \co
\pi_1(M_0,x_0) \to \pi$ and define $f_k=f \circ (h_k \circ \cdots \circ h_1)_*^{-1}\co  \pi_1(M_k,x_k) \to \pi$ for any
$1 \leq k \leq n$. Since $f=h_* \circ f'$, we have that $f_n=f'$.

We arbitrarily orient the upper circles $u^k_i$ and the lower circles $l^k_i$ (so that each Heegaard diagram $D_k$ is
oriented) and denote by $\al^k=(\al^k_1, \cdots, \al^k_{g_k})$ the coloration of the diagram  $D_k$ by the homomorphism
$f_k$.

Up to applying some moves of type II or to well-choosing the orientation of the added circles in a move of type C (or its
inverse), we can assume that the orientation of the upper and lower circles are transported by the homeomorphisms $h_i$.
Note that if we change the orientation of an upper circle $u^k_i$ to its inverse, then the color $\al^k_i=f([\ga^k_i])$
is replaced by $f([(\ga^k_i)^{-1}])=(\al^k_i)^{-1}$, where $\ga^k_i$ is a loop on $(S_k,x_k)$ which crosses (in a
positively-oriented way) the upper circle $u^k_i$ in exactly one point and does not intersect any other upper circle.

We have to verify that, for any $1 \leq k \leq n$, the colors of the diagrams $D_{k-1}$ and $D_k$ are related as
described in the moves of type I-V of Section~\ref{s:colheegdiags}. Without loss of generality, we can assume that $n=1$,
that is that the diagrams $D=(S,u,l)$ colored by $f$ and $D'=(S',u',l')$ colored by $f'$ are related by moves of type
I-V. Denote the genus and color of $D$ (resp.\@ $D'$) by $g$ (resp.\@ $g'$) and $\al=(\al_1, \cdots,\al_g)$ (resp.\@
$\al'=(\al'_1, \cdots,\al'_{g'})$) respectively.

Suppose that $D'$ is obtained from $D$ by a move of type A. Let $1 \leq i \leq g'=g$ and $\ga_i$ be a loop on $(S,x)$
which crosses (in a positively-oriented way) the upper circle $u_i$ in exactly one point and does not intersect any other
upper circle of $D$. Then $\ga'_i=h(\ga_i)$ is a loop on $(S',x')$ which crosses (in a positively-oriented way) the upper
circle $h(u_i)=u'_i$ in exactly one point and does not intersect any other upper circle of $D'$. Therefore
\begin{equation*}
  \al'_i = f'([\ga'_i]) = f'([h ( \ga_i)]) =  f' \circ (h)_* ([\ga_i]) =  f([\ga_i])
   = \al_i.
\end{equation*}
Hence the \p colored Heegaard diagrams $D$ and $D'$ are related by a move of type~I.

Suppose that $D'$ is obtained from $D$ by a move of type B. Then the colors of the upper circles $u'_i$ and $u_i$ agree
(by the same argument as above, since $S'=h(S)$ and $u'=h(u)$). Therefore the \p colored Heegaard diagrams $D'$ is
obtained from the \p colored Heegaard diagram $D$ by a finite sequence of move of type I and III (by decomposing the
isotopy into two-point moves, see Section~\ref{s:colheegdiags}).

Suppose that $D'$ is obtained from $D$ by a move of type C. Since $u'=h(u) \cup \{C_1\}$, the color of the upper circle
$u_i'=h(u_i)$ agrees with that of the upper circle $u_i$ for any $1 \leq i \leq g=g' -1$. Let $\ell$ be a path connecting
the point $x'$ to the circle $C_2$ which does not intersect any upper circle of $D'$. Then the loop $\ell C_2 \ell^{-1}$
crosses $C_1$ in exactly one point and does not intersect any other upper circle of $D'$. Set $\nu=+1$ if $\ell C_2
\ell^{-1}$ crosses $C_1$ in a positively-oriented way and $\nu=-1$ otherwise. Therefore
\begin{equation*}
 \al'_{g'}=f'([\ell C_2^\nu \ell^{-1}])=f'([\ell C_2 \ell^{-1}])^\nu.
\end{equation*}
Now the circle $C_2$ is contractible in $M'$. Thus $[\ell C_2 \ell^{-1}]=1 \in \pi_1(M',x')$ and so $\al'_{g'}=1 \in
\pi$. Hence the \p colored Heegaard diagram $D'$ is obtained from the \p colored Heegaard diagram $D$ by a move of type I
and then a move of type IV.

Suppose that $D'$ is obtained from $D$ by a move of type D. Since $S'=h(S)$ and $u'=h(u)$, the colors of the upper
circles of $D'$ and $D$ agree. Then the \p colored Heegaard diagram $D'$ is obtained from the \p colored Heegaard diagram
$D$ by a move of type I and then a move of type V.

Finally, suppose that $D'$ is obtained from $D$ by a move of type E, i.e., that $u'$ is obtained from $h(u)$ by sliding a
circle $h(u_i)$ past another circle $h(u_j)$. Let $b \co  I \times I \to S'$ be a band which connects $h(u_i)$ to
$h(u_j)$ (that is, such that $b(I \times I) \cap h(u_i)=b(0 \times I)$ and $b(I \times I) \cap h(u_j)=b(1 \times I)$) but
does not intersect any other circle. We can also assume that $x' \not \in b(I \times I)$. Then
\begin{equation*}
u'_i=h(u_i) \#_b h(u_j)=h(u_i) \cup h(u_j) \cup b(I \times \partial I) \setminus
  b(\partial I \times I)
\end{equation*}
and $u'_j$ is a copy $h(u_j)$ which is slightly isotoped such that it has no point in common with $u'_i$. Set
$p=b(0,\frac{1}{2}) \in h(u_i)$ and $q=b(1,\frac{1}{2}) \in h(u_j)$. Up to first applying a move of type II to $u_i$
and/or $u_j$, we can assume that the basis $(d_p b( \cdot,\frac{1}{2}), d_p h(u_i))$ for $T_pS'$ is negatively-oriented
and that the basis $(d_q b( \cdot,\frac{1}{2}), d_q h(u_j))$ for $T_qS'$ is positively-oriented. Then the orientations of
$u'_i$ induced by $h(u_i)$ and $h(u_j)$ coincide and $u'_i$ is endowed with this orientation. Let $\ga_i$ (resp.\@
$\ga_j$) be a loop on $(S,x)$ which crosses (in a positively-oriented way) the upper circle $u_i$ (resp.\@ $u_j$) in
exactly one point and does not intersect any other upper circle of $D$ neither the band $h^{-1}(b(I \times I))$. Let
$\ell_1\co I \to S'$ be a path with $\ell_1(0)=x'$ and $\ell_1(1)=p$ which does not intersect any upper circle of $D'$
and such that $(d_p \ell_1, d_p h(u_i))$ is a negatively-oriented basis for $T_pS'$. Let $\ell_2\co I \to S'$ be a path
with $\ell_2(0)=q$ and $\ell_2(1)=x'$ which does not intersect any upper circle of $D'$ and such that $(d_q \ell_2, d_q
h(u_j))$ is a positively-oriented basis for $T_qS'$, see Figure~\ref{loopab}.
\begin{figure}[t]
                         \psfrag{a}[][]{$\ell_1$}
                         \psfrag{b}[][]{$\ell_2$}
                         \psfrag{q}[][]{$q$}
                         \psfrag{p}[][]{$p$}
                         \psfrag{g}[][]{$h(\ga_i)$}
                         \psfrag{j}[][]{$h(\ga_j)$}
                         \psfrag{t}[][]{$b(\cdot,\frac{1}{2})$}
                         \psfrag{u}[][]{$u'_i$}
                         \psfrag{v}[][]{$u'_j$}
                         \psfrag{x}[][]{$x'$}
                         \scalebox{.9}{\includegraphics{loopab.eps}}
     \caption{}
     \label{loopab}
\end{figure}

Set $\ga'_i=h(\ga_i)$ (resp.\@ $\ga'_j=\ell_2 b(\cdot,\frac{1}{2}) \ell_1$). It is a loop on $(S',x')$ which crosses (in
a positively-oriented way) the upper circle $u'_i$ (resp.\@ $u'_j$) in exactly one point and does not intersect any other
upper circle of $D'$. Therefore we have
\begin{equation*}
  \al'_i = f'([\ga'_i]) = f'([h ( \ga_i)]) =  f' \circ (h)_* ([\ga_i]) =  f([\ga_i])
   = \al_i
\end{equation*}
and, since $\ga'_j$ is homotopic (in $M'$) to the loop $ h(\ga_i)^{-1}h(\ga_j)$,
\begin{eqnarray*}
  \al'_j&  = & f'([\ga'_j])\\
         &  = & f'([ h(\ga_i)^{-1}h(\ga_j)]) \\
         &  = & f' \circ (h)_* ([ (\ga_i)^{-1}\ga_j])\\
         &  = & f([\ga_i]^{-1}[\ga_j]) \\
         &  = & (\al_i)^{-1} \al_j.
\end{eqnarray*}
Hence the \p colored Heegaard diagram $D'$ is obtained from the \p colored Heegaard diagram $D$ by a move of type
I and then a move of type V.\\

Let us prove Part (b.i) of Theorem~\ref{invKuth1}. Let $\xi=(p\co \Tilde{M} \to M)$ be a flat \p bundle  over a \trois
manifold and $\tilde{x}$, $\tilde{x}'$ be two points in $\Tilde{M}$ which belong to the same path-component. Consider a
path $\tilde{\ga}$ in $\Tilde{M}$ connecting $\tilde{x}$ to $\tilde{x}'$. Pushing $\tilde{x}$ to $\tilde{x}'$ along
$\tilde{\ga}$ inside a tubular neighborhood of $\mathrm{Im}(\tilde{\ga})$ in $\Tilde{M}$ yields a self-homeomorphism of
$\Tilde{M}$ which is an equivalence between the pointed flat \p bundles $(\xi,\tilde{x})$ and
$(\xi,\tilde{x}')$. Therefore $K_H(\xi,\tilde{x})=K_H(\xi,\tilde{x}')$ by Part (a) of Theorem~\ref{invKuth1}.\\

Let us prove Part (b.ii) of Theorem~\ref{invKuth1}. Fix a flat \p bundle $\xi=(p\co \Tilde{M} \to M)$ over a \trois
manifold $M$. Let $\tilde{x}$,  $\tilde{x}'$ be two points in $\Tilde{M}$. Set $x=p(\tilde{x})$ and $x'=p(\tilde{x}')$.
Since $M$ is connected, there exists a path $\ga$ in $M$ connecting $x$ to $x'$. Let $\tilde{z}\in \Tilde{M}$ be the
endpoint of the lift of $\ga$ to a path that starts at $\tilde{x}$. Since $p(\tilde{z})=x'=p(\tilde{x}')$, there exists
$\be \in \pi$ such that $\tilde{x}'=\be \cdot \tilde{z}$.

Firstly, suppose that the monodromy $f$ of $\xi$ at $\tilde{z}$ is surjective. Then there exists a loop $\sigma$ based on
$x'$ such that $f([\sigma])=\be$. Denote by $\tilde{\sigma}$ the lift of $\sigma$ to a path that starts at $\tilde{z}$.
Since $\tilde{x}'=\be\cdot \tilde{z}= f([\sigma]) \cdot \tilde{z}$, the path $\tilde{\sigma}$ ends at $\tilde{x}'$.
Finally, $\tilde{\ga} \tilde{\sigma}$ is a path in $\Tilde{M}$ which connects $\tilde{x}$ to $\tilde{x}'$. Hence
$K_H(\xi,\tilde{x})=K_H(\xi,\tilde{x}')$ by Part (b.i) of Theorem~\ref{invKuth1}.

Secondly, suppose that the Hopf \p coalgebra $H$ is crossed (see Section~\ref{deficro}). Since $\tilde{x}'=\be \cdot
\tilde{z}$, the monodromies $f_{\tilde{z}},\, f_{\tilde{x}'}\co \pi_1(M,x') \to \pi$ of $\xi$ at $\tilde{z}$ and
$\tilde{x}'$ are related by $f_{\tilde{x}'}=\be f_{\tilde{z}}\, \bei$. Let $D=(S,u,l)$ be a Heegaard diagram of genus $g$
of $M$ whose upper and lower circles are arbitrarily oriented. Denote by $D_{f_{\tilde{z}}}$ and $D_{f_{\tilde{x}'}}$ the
\p colored Heegaard diagrams obtained by coloring $D$ with $f_{\tilde{z}}$ and $f_{\tilde{x}'}$. Since
$f_{\tilde{x}'}=\be f_{\tilde{z}}\, \bei$, the colors of $D_{f_{\tilde{z}}}$ and $D_{f_{\tilde{x}'}}$ are conjugate.
Therefore:
\begin{eqnarray*}
  K_H(\xi,\tilde{x}') & = & K_H(D_{f_{\tilde{x}'}}) \\
               & = & K_H(D_{f_{\tilde{z}}}) \text{ \quad by Lemma~\ref{cojugclassindepKH}} \\
               & = & K_H(\xi,\tilde{z})  \\
               & = & K_H(\xi,\tilde{x}) \text{ \quad by Part (b.i) of Theorem~\ref{invKuth1}.}
\end{eqnarray*}

Finally, suppose that $\pi$ is abelian. Then the Hopf \p coalgebra $H=\{H_\al\}_{\al \in \pi}$ is crossed (see
Section~\ref{deficro}) and so $K_H(\xi,\tilde{x})=K_H(\xi,\tilde{x}')$ by the previous case. This completes the proof of
Theorem~\ref{invKuth1}.

\section{Examples}\label{sect-kup3}

In this section, we give some examples of computations of the scalar invariant of flat bundle over \trois manifolds
constructed in Section~\ref{sect-kup2}.

\subsection{Example}
As remarked by Vainerman~\cite{Vai1}, the Kac-Paljutkin Hopf algebra $A=\mathbb{C}^4 \oplus \mathrm{Mat}_2(\mathbb{C})$,
viewed as a central prolongation of $F(\mathbb{Z}/2\mathbb{Z})$, leads to a finite type involutory Hopf
$\mathbb{Z}/2\mathbb{Z}$-coalgebra $H=\{H_0,H_1\}$ over $\mathbb{C}$. Namely, set $H_0= \mathbb{C}^4$ and
$H_1=\mathrm{Mat}_2(\mathbb{C})$ as algebras. Let $\{e_1, e_2, e_3, e_4 \}$ be the (standard) basis of $H_0$ and
$\{e_{1,1},e_{1,2},e_{2,1},e_{2,2}\}$ be the (standard) basis of $H_1$. The counit $\varepsilon\co H_0 \to \mathbb{C}$ is
given by $\varepsilon(e_1)=1$ and $\varepsilon(e_2)=\varepsilon(e_3)=\varepsilon(e_4)=0$. The comultiplication is given
by
\begin{equation*}
\begin{array}{l}
   \cp{0}{0}(e_1)=e_1 \otimes e_1 + e_2 \otimes e_2 + e_3 \otimes e_3 + e_4 \otimes e_4  \\
   \cp{0}{0}(e_2)=e_1 \otimes e_2 + e_2 \otimes e_1 + e_3 \otimes e_4 + e_4 \otimes e_3  \\
   \cp{0}{0}(e_3)=e_1 \otimes e_3 + e_3 \otimes e_1 + e_2 \otimes e_4 + e_4 \otimes e_2  \\
   \cp{0}{0}(e_4)=e_1 \otimes e_4 + e_4 \otimes e_1 + e_2 \otimes e_3 + e_3 \otimes e_2
\end{array}
\end{equation*}
\begin{equation*}
\begin{array}{l}
   \cp{0}{1}(e_{1,1})
   =e_1\otimes e_{1,1} + e_2\otimes e_{2,2} + e_3 \otimes e_{1,1} + e_4 \otimes e_{2,2} \\
   \cp{0}{1}(e_{1,2})
   =e_1 \otimes e_{1,2} -i\, e_2 \otimes e_{2,1} - e_3 \otimes e_{1,2} +i\, e_4 \otimes e_{2,1} \\
   \cp{0}{1}(e_{2,1})
   =e_1 \otimes e_{2,1} +i\, e_2 \otimes e_{1,2} - e_3 \otimes e_{2,1} -i\, e_4 \otimes e_{1,2} \\
   \cp{0}{1}(e_{2,2})
   =e_1 \otimes e_{2,2} + e_2 \otimes e_{1,1} + e_3 \otimes e_{2,2} + e_4 \otimes e_{1,1}
\end{array}
\end{equation*}
\begin{equation*}
\begin{array}{l}
   \cp{1}{0}(e_{1,1})
   =e_{1,1}\otimes e_1 + e_{2,2}\otimes e_2 + e_{1,1} \otimes e_3 + e_{2,2} \otimes e_4 \\
   \cp{1}{0}(e_{1,2})
   =e_{1,2}\otimes e_1 +i\, e_{2,1}\otimes e_2 - e_{1,2} \otimes e_3 -i\, e_{2,1} \otimes e_4 \\
   \cp{1}{0}(e_{2,1})
   =e_{2,1}\otimes e_1 -i\, e_{1,2}\otimes e_2 - e_{2,1} \otimes e_3 +i\, e_{1,2} \otimes e_4 \\
   \cp{1}{0}(e_{2,2})
   =e_{2,2}\otimes e_1 + e_{1,1}\otimes e_2 + e_{2,2} \otimes e_3 + e_{1,1} \otimes e_4
\end{array}
\end{equation*}
\begin{equation*}
\begin{array}{l}
   \cp{1}{1}(e_1)
   =\frac{1}{2}\,(e_{1,1} \otimes e_{1,1} + e_{2,2} \otimes e_{2,2} + e_{1,2} \otimes e_{1,2} + e_{2,1} \otimes e_{2,1}) \\
   \cp{1}{1}(e_2)
   =\frac{1}{2}\,(e_{1,1} \otimes e_{2,2} + e_{2,2} \otimes e_{1,1} + i\, e_{1,2} \otimes e_{2,1} -i\, e_{2,1} \otimes e_{1,2})\\
   \cp{1}{1}(e_3)
   =\frac{1}{2}\,(e_{1,1} \otimes e_{1,1} + e_{2,2} \otimes e_{2,2} - e_{1,2} \otimes e_{1,2} - e_{2,1} \otimes e_{2,1}) \\
   \cp{1}{1}(e_4)
   =\frac{1}{2}\,(e_{1,1} \otimes e_{2,2} + e_{2,2} \otimes e_{1,1} -i\, e_{1,2} \otimes e_{2,1} +i\, e_{2,1} \otimes e_{1,2})
\end{array}
\end{equation*}
The antipode is given by $S_0(e_k)=e_k$ for any $1 \leq k \leq 4$ and $S_1(e_{k,l})=e_{l,k}$ for any $1 \leq k,l \leq 2$.

Let $n\geq 1$. There exists two flat $\mathbb{Z}/2\mathbb{Z}$-bundles $\xi_n^0$ and $\xi_n^1$ over the lens space
$L(2n,1)$, whose monodromies $f_n^0, f_n^1\co  \pi_1(L(2n,1)) \cong \mathbb{Z}/2n\mathbb{Z} \to \mathbb{Z}/2\mathbb{Z}$
are respectively given by $f_n^0(1)=0$ and $f_n^1(1)=1$.

A Heegaard diagram $\{u_1,l_1\}$ of genus $1$ of the lens space $L(2n,1)$ is given, on the torus
$\mathbb{T}=\mathbb{R}^2/\mathbb{Z}^2$, by $u_1=\mathbb{R}(0,1) + \mathbb{Z}^2$ and $l_1=\mathbb{R}(1,\frac{1}{2n}) +
\mathbb{Z}^2$. Fix $k =0,1$ and set $\al=f_n^k(1)\in \mathbb{Z}/2\mathbb{Z}$. Denote by $D_\al$ the \p colored Heegaard
diagram obtained from $(\mathbb{T},\{u_1,l_1\})$ by providing the circle $u_1$ with the color $\al$. Then
$$K_H(\xi_n^k)=(\dim H_0)^{-1} K_H(D_\al)=\frac{1}{4} K_H(D_\al),$$ where $K_H(D_\al) \in \mathbb{C}$ equals the tensor
depicted in Figure~\ref{lens1}.
\begin{figure}[h]
     \subfigure[$K_H(D_\al)$]{\label{lens1}
                               \psfrag{d}[Br][Br]{$\Delta_{\underbrace{\al, \dots, \al}_{2n \mathrm{ \; times}}}$}
                               \psfrag{m}{$m_\al$}
                               \scalebox{.9}{\includegraphics{lensex1b.eps}}} \hspace*{1.3cm}
     \subfigure[$F_\al\co H_0 \to H_\al$]{\label{lens2}
                               \psfrag{d}[Br][Br]{$\Delta_{\underbrace{\al, \dots, \al}_{8 \mathrm{ \; times}}}$}
                               \psfrag{m}{$m_\al$}
                               \scalebox{.9}{\includegraphics{lensex2b.eps}}}
     \caption{}
     \label{lens0}
\end{figure}

Let $F_\al\co H_0 \to H_\al$ be the map defined  in Figure~\ref{lens2}. One easily checks by hand that $F_\al(x)=
\varepsilon(x) \, 1_\al$ for all $x \in H_0$. Then, using the (co)associativity of the (co)multiplication, we get the
equalities of Figure~\ref{lensfin}.
\begin{figure}[h]
                               \psfrag{d}[Br][Br]{$\Delta_{\underbrace{\al, \dots, \al}_{2(n+4) \mathrm{ \; times}}}$}
                               \psfrag{m}{$m_\al$}
                               \psfrag{b}[Br][Br]{$\Delta_{\underbrace{\al, \dots, \al}_{2n \mathrm{ \; times}}}$}
                               \psfrag{s}[Br][Br]{$\Delta_{0,\underbrace{\al, \dots, \al}_{2n \mathrm{ \; times}}}$}
                               \psfrag{u}[Br][Br]{$\Delta_{\underbrace{\al, \dots, \al}_{8 \mathrm{ \; times}}}$}
                               \psfrag{e}{$\varepsilon$}
                               \psfrag{1}[Br][Br]{$1_\al$}
                               \psfrag{=}{$=$}
                               \scalebox{.9}{\includegraphics{lensex3a.eps}}
     \caption{}
     \label{lensfin}
\end{figure}

Therefore $K_H(\xi_{n+4}^k)= K_H(\xi_n^k)$ for any $n \geq 1$. Hence, by computing by hand the values of $K_H(\xi_n^k)$
for $1\leq n \leq 4$ and $0 \leq k \leq 1$, we obtain that :
\begin{equation*}
K_H(\xi_n^0)=4 \text{\quad and \quad} K_H(\xi_n^1)=\left \{ \begin{array}{ll} 2 & \text{ if $n \equiv 1 \mathrm{\; (mod
\;} 2)$}\\ 4 & \text{ if $n \equiv 0 \mathrm{\; (mod \;} 4)$}\\ 0 & \text{ if $n \equiv 3 \mathrm{\; (mod \;} 4)$}
\end{array}\right. .
\end{equation*}

\subsection{Example}
Let $\pi$ and $G$ be two finite groups, and $\phi\co G \to \pi$ be a group homomorphism. Then $\phi$ induces a Hopf
algebras morphism $F(\pi) \to F(G)$, given by $f \mapsto f \circ \phi$, whose image is central. Here $F(G)=\mathbb{C}^G$
and $F(\pi)=\mathbb{C}^\pi$ denote the Hopf algebras of complex-valued functions on $G$ and $\pi$ respectively. By
Section~\ref{abstrasctreformulation}, this data yields to a Hopf \p coalgebra $H^\phi=\{H^\phi_\al\}_{\al \in \pi}$.

Let us describe more precisely this Hopf \p coalgebra. Denote by $(e_g)_{g \in G}$ the standard basis of $F(G)$ given by
$e_g(h)=\delta_{g,h}$. Then, for any $\al,\be \in \pi$, we have that:
\begin{equation*}
\begin{array}{l}
\displaystyle  H^\phi_\al  = \!\!\!\!\!\sum_{g \in \phi^{-1}(\al)}\!\!\!\! \mathbb{C}\, e_g, \qquad
  m_\al(e_g \otimes e_h)  = \delta_{g,h} \, e_g \text{ \, for any $g,h \in \phi^{-1}(\al)$,} \\[15pt]
\displaystyle  1_\al  = \!\!\!\!\!\sum_{g \in \phi^{-1}(\al)} \!\!\!\! e_g, \qquad \quad \,
\varepsilon(e_g)  =\delta_{g,1} \text{ \, for any $g\in \phi^{-1}(1)$,} \\[15pt]
\displaystyle  \cp{\al}{\be}(e_g) = \!\!\!\!\!  \sum_{\shortstack[c]{ $\scriptstyle h\in \phi^{-1}(\al)$ \\ $\scriptstyle
k \in \phi^{-1}(\be)$  \\ $\scriptstyle hk=g$}} \!\!\!\! e_h \otimes e_k
  \text{ \, for any $g \in \phi^{-1}(\al\be)$,} \\[35pt]
\displaystyle S_\al(e_g) =e_{g^{-1}} \text{ \, for any $g\in \phi^{-1}(\al)$.}\\[5pt]
\end{array}
\end{equation*}

Note that the Hopf \p coalgebra $H^\phi=\{H^\phi_\al\}_{\al \in \pi}$ is involutory and of finite type. Since $\dim
H^\phi_1=\# \phi^{-1}(1_G)=\# \ker \phi$ is non-zero in the field $\mathbb{C}$, the invariant $K_{H^\phi}$ of pointed
flat \p bundles over \trois manifolds is well defined.
\begin{lemma}
Let $(\xi,\tilde{x})$ be a pointed flat \p bundle over a \trois manifold $M$. Then $$K_{H^\phi}(\xi,\tilde{x})=\#\{g\co
\pi_1(M,x) \to G \;|\; \phi \circ g=f \},$$ where $x$ is the image of $\tilde{x}$ under the bundle map $\Tilde{M} \to M$
of $\xi$ and $f\co \pi_1(M,x) \to \pi$ is the monodromy of $\xi$ at~$\tilde{x}$.
\end{lemma}
\begin{proof}
Using the above explicit description of $H^\phi$, one easily gets that, for any $\al \in \pi$ and $g_1, \dots, g_m \in
\phi^{-1}(\al)$,
\begin{equation*}
 \psfrag{m}{$m_\al=\begin{cases} 1 & \text{if } g_1= g_2= \cdots = g_m \\ 0& \text{otherwise} \end{cases},$}
 \psfrag{c}[][]{$e_{g_1}$}
 \psfrag{a}[][]{$e_{g_m}$}
 \scalebox{.9}{\includegraphics{exp1.eps}}
\end{equation*}
and that, for any $\al_1, \dots, \al_n \in \pi$ such that $\al_1 \cdots \al_n=1$ and $g_1 \in \phi^{-1}(\al_1), \dots,
g_n \in \phi^{-1}(\al_n)$,
\begin{equation*}
 \psfrag{u}[][]{$e_{g_n}$}
 \psfrag{c}[][]{$e_{g_1}$}
 \psfrag{d}[][]{$\Delta_{\al_1, \dots, \al_n}$}
 \psfrag{a}{$= \,(\dim H^\phi_1) \;\delta_{g_1 \cdots g_n, 1}$.}
 \scalebox{.9}{\includegraphics{exp2.eps}}
\end{equation*}

Let $(\xi,\tilde{x})$ be a pointed flat \p bundle over a \trois manifold $M$. Let $D=(S,u,l)$ be a Heegaard diagram of
genus $g$ of $M$. Orient it and color it by the monodromy $f$ of $\xi$ at $\tilde{x}$. Denote by $\al=(\al_1,
\dots,\al_g)$ the color of $D$. Let $w_1, \dots, w_g$ be the words in the alphabet $\{X_1, \dots, X_g\}$ as in
Section~\ref{s:colheegdiags}. For any $1\leq i \leq g$, write $ w_i(X_1, \dots, X_g)=X^{\epsilon_{i,1}}_{k_{i,1}} \cdots
X^{\epsilon_{i,n_i}}_{k_{i,n_i}} $ where $\epsilon_{i,j} =\pm 1$ and $1 \leq k_{i,j} \leq g$.

The non-zero factors of the tensor associated to an upper circle $u_i$ are of the form:
\begin{equation*}
 \psfrag{m}{$m_{\al_i}=1$ \; where \; $g_i \in \phi^{-1}(\al_i)$.}
 \psfrag{c}[][]{$e_{g_i}$}
 \psfrag{a}[][]{$e_{g_i}$}
 \scalebox{.9}{\includegraphics{exp4.eps}}
\end{equation*}
The non-zero factors of the tensor associated to a lower circle $l_i$ are of the form:
\begin{equation*}
 \psfrag{u}{$e_{h^i_1}$}
 \psfrag{c}{$e_{h^i_{n_i}}$}
 \psfrag{d}[][]{$\Delta_{\al^{\epsilon_{i,1}}_{k_{i,1}}, \dots, \al^{\epsilon_{i,n_i}}_{k_{i,n_i}}}$}
 \psfrag{a}{$=\, (\dim H^\phi_1)\; \delta_{h^i_1 \cdots h^i_{n_i},1}$ \;  where
 \; $h^i_j \in \phi^{-1}(\al^{\epsilon_{i,j}}_{k_{i,j}})$.}
 \scalebox{.9}{\includegraphics{exp6.eps}}
\end{equation*}

Now, at each crossing, that is, for $1 \leq i \leq g$ and $1 \leq j \leq n_i$, the contraction rule amounts to the
equality of $h^i_j$ with $g^{\epsilon_{i,j}}_{k_{i,j}}$. Therefore we get that:

\begin{eqnarray*}
\lefteqn{K_{H^\phi}(\xi,\tilde{x})}\\
  & = & (\dim H^\phi_1)^{-g} \!\!\!\!\!\! \sum_{\shortstack[c]{ $\scriptstyle g_i\in \phi^{-1}(\al_i)$ \\ $\scriptstyle h^i_j
        \in \phi^{-1}(\al^{\epsilon_{i,j}}_{k_{i,j}})$}} \!\!\! \Bigl (\prod_{1\leq i \leq g} (\dim H^\phi_1)\;
        \delta_{h^i_1 \cdots
        h^i_{n_i},1}   \Bigr ) \Bigl ( \prod_{\shortstack[c]{$\scriptstyle 1 \leq i \leq g$ \\ $\scriptstyle 1 \leq j \leq n_i$}}
        \delta_{h^i_j, g^{\epsilon_{i,j}}_{k_{i,j}}} \Bigr ) \\
  & = & \sum_{g_i\in \phi^{-1}(\al_i)} \Bigl (  \prod_{1\leq i \leq g} \delta_{w_i(g_1, \dots, g_g),1}   \Bigr )\\
  & = & \#\{(g_1, \dots, g_g) \in G^g \;|\; w_i(g_1, \dots, g_g)=1 \text{ and } \phi(g_i)=\al_i \; \forall \, 1 \leq i \leq g\} \\
  & = & \#\{g\co \pi_1(M,x) \to G \;|\; \phi \circ g=f \}.
\end{eqnarray*}
\end{proof}

For example, suppose that $\pi=\mathbb{Z}/2\mathbb{Z}$, $G=S_3$, and $\phi\co S_3 \to \mathbb{Z}/2\mathbb{Z}$ is the
signature of \trois permutations. Let $M$ be a closed, connected, and oriented \trois manifold such that
$\pi_1(M)=\mathbb{Z}/2\mathbb{Z}$. There is two flat $\mathbb{Z}/2\mathbb{Z}$-bundles $\xi_0$ and $\xi_1$ over $M$, whose
monodromies are respectively trivial and $\id_{\mathbb{Z}/2\mathbb{Z}}$. In this setting, we get that
$K_{H^\phi}(\xi_0)=1$ and $K_{H^\phi}(\xi_1)=3$.

\bibliographystyle{amsplain}
\bibliography{kup}
\end{document}